\documentclass[11pt]{article}
\usepackage{fancybox,%
            fancyheadings,%
            epsf,%
            theorem,%
            eepic,%
            psfig,%
            epsfig,%
            amssymb,%
            amsmath,%
            pifont,%
            subeqn
}

\usepackage%
[
nodraft,
]
{symbols}
 
 \let\c=\cedille
                  
 \def\coauthoronesname{D}
 
 \def\coauthortwoname{A}

\renewcommand\norm[1]{\left|#1\right|}

\newtheorem{stmt}{\hspace{-1mm}(\hspace{-1mm}}

\newenvironment{statement}{\pullup{7mm}\begin{description}\item[]\begin{stmt}{\em \shiftleft{2mm})}}{\end{stmt}\end{description}}

\renewcommand{\Re}{\mathbb R}

\def\BD{{\cal B}(\Delta)}

\def\Bmu{{\cal B}(\mu)}
\def\icxbr{(\t\circ,\x\circ)\in\mR\times B_r}

\begin{document}
\title{Persistency of excitation for uniform convergence \\ in nonlinear control systems\footnote{This work is supported in part by a CNRS-NSF collaboration project and was partially done while the first two authors were visiting the CCEC at UCSB. It is also supported by the AFOSR
under grant F49620-00-1-0106 and the NSF under grants ECS-9988813 and INT-9910030.}}
\author{Antonio Lor\'{\i}a$^\dagger$\quad %
Elena Panteley$^*$ \quad %
Dobrivoje Popovi\'c$^\bullet$ \quad %
Andrew R. Teel$^\bullet$ \quad %
\address{$^\dagger$C.N.R.S, UMR 5228, Laboratoire d'Automatique de Grenoble,
ENSIEG, St. Martin d'H\`eres, France
\\[2mm]
$^*$ IPME, Acad. of Sciences of Russia
 61, Bolshoi Ave,
 St. Petersburg, Russia
 elena@ieee.org
\\[2mm]
$^\bullet$Center for Control Engineering and Computation (CCEC),  Department of Electrical and Computer Engineering, University  of California,  Santa  Barbara, CA 93106-9560, USA
}}
\date{  %
 { \small \today }
 { \small \em -- Preprint --  }
 }
\maketitle
\sloppy

\begin{abstract} 
\small In previous papers we have introduced a sufficient condition for uniform attractivity of the origin for a class of nonlinear time-varying systems which is stated in terms  of {\em persistency of excitation} (PE), a concept well known in the adaptive control and systems identification literature. The novelty of our condition, called uniform $\delta$-PE, is that it is tailored  for nonlinear functions of time and {\em state} and it allows us to prove {\em uniform} asymptotic stability. In this paper we present a new definition of u\ped\ which is conceptually similar to but technically different from its predecessors and give several useful characterizations. We make connections between this property and similar properties previously used in the literature.  We also show when this condition  is necessary and
sufficient for uniform (global)  asymptotic stability for a large class of nonlinear time-varying systems. Finally, we show the utility of our main results on some control applications regarding feedforward systems and systems with matching nonlinearities.
\end{abstract}


\section{Introduction}

In many interesting nonlinear control problems, the closed-loop
control system can be modeled by the differential equation
\begin{equation}
\label{one}
\dot x= F(t,x)
\end{equation}
where $F(\cdot,x)$ is not necessarily periodic.

When convergence of the trajectories of (\ref{one}) to
a given fixed point is required (notice that this includes many problems of trajectory tracking control if we think of $x$ as the tracking error), perhaps the most appealing notion that
includes such convergence is {\em uniform asymptotic stability} (UAS) because
of its inherent robustness and its assertion that the convergence
rate does not depend on the initial time in any significant way.
Such properties are typically not guaranteed when using
Barbalat's lemma, for example, which is often used in adaptive
and nonlinear control to establish convergence.

Nevertheless, one can precisely characterize the  property of UAS in many different ways. For instance, via so-called ${\cal K}{\cal L}$ estimates that is, using a bounding function over the norm of the solutions, which decreases uniformly with time and increases uniformly with the size of the initial states. Another common characterization is via the use of  smooth positive definite Lyapunov functions with a negative definite total derivative. See for instance \cite{HALE,KHALIL}. While these characterizations are very useful as intermediate steps in proving other properties (as for instance robust stability of a perturbed system) in general, Lyapunov  functions or ${\cal K}{\cal L}$ estimates are   difficult to obtain. 

Other sufficient and necessary conditions for uniform asymptotic stability involve verifying certain {\em integral} criteria which are somehow in the spirit of Barbalat's lemma but allow to conclude uniform convergence. For instance, roughly speaking  one can prove that  if the norm of the solution $x(\cdot)$ of \rref{one} and  the integral of this norm squared are bounded by an increasing function of the initial states, $\alpha(\norm{x(t_\circ)})$, which is zero only at zero, then the solutions necessarily converge to zero. One may also formulate conditions using more generic functions (other than quadratic) of the norm of the solutions. Furthermore, such conditions may be stated in a fairly general manner to cover stability of sets (not necessarily compact) other than the origin. Such sufficient and necessary conditions have been presented recently in \cite{INTLEMMCSS} for systems described by differential inclusions. See also the references therein. 

In order to meet the required conditions imposed on the integral equations mentioned above or  on the derivatives of Lyapunov functions, other approaches use observability-type arguments and $exciteness$ conditions. Roughly speaking one tries to identify a converging output and then verifies whether all the modes of the zero-dynamics (i.e., the dynamics which is left by zeroing the output) are sufficiently excited so that all converge to zero. For linear systems such methods have been under investigation, based on \cite{KALMAN}, starting probably with \cite{ASTBOH} and followed by numerous  works including \cite{AND77,MORNAR2,SASBOD,NARANA}.  In the nonlinear case we find for instance \cite{AEYSEPPEUT,ARST82,TACDELTAPE}. See also the interesting recent paper \cite{TCLEETAC01} which establishes sufficient conditions for uniform convergence in terms of a detectability notion and appealing to the technique of limiting equations\footnote{Roughly speaking, dynamic equations describing the limiting (as $t\to\infty$) behaviour of the system.}.

The advantage of methods using observability-type arguments and related $exciteness$ conditions is that for a certain class of systems one may infer the stability and convergence properties simply by ``looking'' at the dynamics of the system. For instance, for linear systems
\[
\dot x = -P(t) x,\qquad  P(t) = P(t)^T \geq 0,\ \forall \, t\geq 0
\]
it was shown in \cite{MORNAR2} that it is necessary and sufficient for uniform asymptotic stability that $P(t)$ be {\em persistently exciting} (PE), namely that there exist $a>0$ and $b\in\mR$ such that for {\em all } unitary vectors $x\in\mR^n$, 
\begin{equation}
\label{pe:mornar}
\int_{t_\circ}^t \norm{P(s)x}ds \geq a(t-t_\circ) + b\qquad \forall\, t\geq t_\circ\geq 0\,. 
\end{equation}
Equivalently, this system is uniformly (in $t_\circ$) completely (i.e. for all initial states) observable (see e.g. \cite{CLAN8}) from the output $y(t) := P(t)x(t) $ if and only if $P(\cdot)$ is PE. Notice that $P(t)$ does not need to be full rank for any fixed $t$. 

That PE is a necessary condition for uniform asymptotic (exponential) stability for linear systems has been well-known for  many years now. In the case of general nonlinear systems, this was established for a generalized notion of persistency of excitation by Artstein in \cite[Theorem 6.2]{ARST78}. Relying on the notion of limiting equations  it was proved in \cite{ARST78} that for the system \rref{one} it is necessary for uniform convergence that for each $\delta>0$ there exist $a>0$ and $b\in\mR$ such that
\begin{equation}
\label{pe:arst}
\norm{x} \geq \delta,\ t\geq t_\circ\geq 0\quad\Rightarrow\quad \int_{t_\circ}^t \norm{F(s,x)}ds \geq a(t-t_\circ) + b\,.
\end{equation}
Sufficiency was also shown under other conditions involving Lyapunov functions and making use of a theorem establishing that uniform asymptotic stability of \rref{one} is equivalent to asymptotic stability of {\em all} the limiting equations of \rref{one}. 

In view of the importance that persistency of excitation and related observability conditions
 have been proved to have in establishing convergence of linear and nonlinear systems, PE has been at the basis of the formulation of  sufficient conditions for uniform asymptotic stability in many contexts.  For examples, see \cite[Theorem 6.3]{ARST78} and the results in \cite{KALMAN},  \cite{ASTBOH}, \cite{AND77,MORNAR2,SASBOD,NARANA}, and \cite{AEYSEPPEUT,ARST82,TACDELTAPE}. In the context of adaptive control of nonlinear parameterized systems, other notions of persistency of excitation for nonlinear systems have been introduced recently in \cite{ANA99}.

In \cite{uped,TACDELTAPE} we introduced a sufficient condition for uniform attractivity for a certain class of nonlinear systems. In words, our condition is that  a certain function $(t,x)\mapsto\phi(t,x)$ evaluated along the trajectories of the system, be persistently exciting whenever the trajectories  are bounded away from a $\delta$-neighbourhood of the origin (cf. Def. \ref{def:olddpe}). Loosely speaking, this property, called uniform $\delta$-persistency of excitation (u\ped), ensures that the zero-dynamics of the system is sufficiently excited. 

In this paper, we present a new definition of u\ped\ which is stated in a form that does not involve the state trajectories (see Def. \ref{def:newdpe}) and therefore it is much simpler to verify. As a matter of fact, it is also a reformulation of the condition \rref{pe:arst} and for the system \rref{one}, we will show both necessity and sufficiency for uniform convergence without relying on limiting equations theory. The former follows simply using a Lipschitz condition and Gronwall's lemma. The proof of sufficiency relies  on a recent version of Matrosov's theorem (see \cite[Theorem 5.5]{ROUMAW}) reported in \cite{MATCDC02}. We will also present some applications in control systems design, including systems with ``matching nonlinearities'' which appear in applications involving passivity-based control, Model Reference Adaptive Control (MRAC), closed-loop identification, feedforward systems, among others. These results, which generalize in certain directions other contributions in the literature, follow as corollaries of our main theorems based on the property of U\ped.

The rest of the paper is organized as follows: in the next section we define our notation and give some preliminary definitions. In Section \ref{sec:defsdpe} we present the definition of uniform $\delta$-PE (u\ped) and related properties. In Section \ref{sec:nec} we give a short proof of necessity of u\ped. In Section \ref{sec:suf} we present a result on sufficiency of u\ped\ which covers Artstein's result. In Section \ref{sec:appl} we present several  application of our main results for a class of systems which includes the closed loop system in MRAC, and closed loop identification of mechanical systems. We conclude with some remarks in Section \ref{sec:conc}.

\section{Preliminaries}

{\bf Notation. } For two constants $\Delta\geq\delta\geq 0$  we define ${\cal H}(\delta,\Delta) \defeq \{x\in \mR^n \, : \, \delta \leq \norm{x} \leq \Delta\}$. We also will use ${\cal B}(r):={\cal H}(0,r)$. A continuous function $\rho:\mRp\to\mRp$ is of class ${\cal N}$ if it is non decreasing. A continuous function $\gamma : \mRp \to \mRp$ is of class ${\cal K}$ ($\gamma
\in {\cal K}$), if it is strictly increasing and $\gamma(0) =
0$; $\gamma \in {\cal K}_\infty$ if in addition, $\gamma(s) \to
\infty$ as $s\to \infty$. A continuous function $\beta : \mRp\times
\mRp \to \mRp$ is of class ${\cal K}{\cal L}$ if $\beta(\cdot,t) \in
{\cal K}$ for each fixed $t\in\mRp$ and $\beta(s,t) \to 0$ as $ t\to
+\infty$ for each $s \geq 0$.  We denote by $x(\cdot,\t\circ,\x\circ)$, the solutions of the differential equation $\dot x = F(t,x)$ with initial conditions $(\t\circ,\x\circ)$. 

We recall that a function $F(\cdot,\cdot)$ is locally Lipschitz in $x$ uniformly in $t$ if for each $x_0$ there exists $L$ such that 
\[
\norm{F(t,x) - F(t,y) } \leq L\norm{x-y}
\]
for all $x$ and $y$ in a neighbourhood of $x_0$ and for all $t\in\mR$. For a locally Lipschitz function $V:\mR\times\mR^n\to \mR$ we define its total time derivative along the trajectories of $\dot x = F(t,x)$ as,  $\dot V(t,x):= \jac{V}{t}+ \jac{V}{x}F(t,x)$\,.  In general, this quantity is defined almost everywhere.

We will consider general nonlinear time-varying systems of the form \rref{one} where $F:\mR \times \mR^n \to \mR^n$ is such that solutions of \rref{one} exist over finite intervals. As it has been motivated e.g. in \cite{TACDELTAPE,MORNAR}, for these systems the most desirable forms of stability are those which are uniform in the initial time:
\begin{defin}[Uniform global stability]
  \label{def:stability}
  The origin of the system \rref{one} is said to be uniformly 
  globally stable (UGS) if there exists $\gamma\in {\cal K}_\infty$ such that, for
  each $(t_{\circ},\x\circ) \in \mR \times \mR^{n}$ each solution $x(\cdot,t_{\circ},x_{\circ})$ satisfies 
\begin{equation}
\label{us}
  \norm{x(t,t_{\circ},x_{\circ})} \leq \gamma(\norm{x_{\circ}})  \qquad \forall \; t \geq t_{\circ} \ .
\end{equation}
\end{defin}
\begin{defin}[Uniform global attractivity]
The origin of the system (\ref{one}) is said to be  uniformly globally attractive if  for each  $r,\,\sigma>0$  there exists $T>0$ such that
\begin{equation}
\label{UA}
\norm{x_{\circ}}\leq r  \, \Longrightarrow \, \|x(t,t_{\circ},x_{\circ})\|
\leq \sigma \qquad  \forall\, t \geq t_{\circ} + T\,
  \ .
\end{equation}
\end{defin}
Furthermore, we say that the (origin of the) system  is uniformly globally asymptotically stable (UGAS) if it is UGS and  uniformly globally  attractive.

We will also make use of the following. 
\begin{definition}[Uniform exponential stability]
The origin of the system (\ref{one}) is said to be  uniformly (locally) 
exponentially stable (ULES) if there exist constants $\ga1, \ga2$ and  $r > 0$  such that for all $\icxbr$ and all corresponding solutions
\begin{equation}
\label{ulesineq}
 \|x(t,t_{\circ},x_{\circ})\| \leq \ga1 \|x_{\circ} \| e^{-\ga2(t-t_{\circ})}
 \qquad\qquad  \forall t\ge\t\circ.
\end{equation}
The system \rref{one} is  uniformly globally exponentially stable (UGES) if there exist $\gamma_1,\;\gamma_2>0$ such that \rref{ulesineq} holds for all $(\t\circ,\x\circ) \in \mR\times \mR^n$.
\end{definition}

\vspace{-.2in}

\section{Uniform $\delta$-persistency of excitation}
\label{sec:defsdpe}
We present a new  definition of ``u\ped'', a property originally introduced in \cite{uped,TACDELTAPE}. The newly defined  property is conceptually similar to each of the previous ones; however, it is technically different in the sense that: first, it is easier to verify since it is formulated as a property inherent to a nonlinear function instead of being directly related to the solutions of a differential equation. Second,  the new property defined below is also necessary for uniform attractivity of the system \rref{one} which is seemingly not the case with the previous definitions. Yet, it is interesting to remark that, in general, neither u\ped\ as defined in \cite{TACDELTAPE} nor as defined below, implies the other.

Let $x\in\mR^n$ be partitioned as $x:=\col[x_1,x_2]$ where $x_1\in\mR^{n_1}$ and $x_2\in\mR^{n_2}$. Define the {\em column vector} function $\phi:\mR\times \mR^n\to \mR^m$ and the set ${\mathcal D}_{1}:= \left(\mR^{n_{1}} \backslash \left\{ 0 \right\}\right) \times \mR^{n_{2}}$.

\vspace{-.1in}

\begin{definition}
\label{def:newdpe}
A function  $\phi(\cdot,\cdot)$ where $t\mapsto\phi(t,x)$ is locally integrable, is said to be uniformly $\delta$-persistently exciting (u\ped) with respect to $x_1$ if 
for each $x \in {\cal D}_1$ there exist $\delta >0$, $T>0$ and $\mu>0$ 
s.t. $\forall t \in \mR$,
\begin{equation}\label{eq:newdpe}
  \norm{z-x} \leq \delta  \quad \Longrightarrow \quad \int_{t}^{t+T} \norm{\phi(\tau,z)} d\tau \geq \mu\,.
\end{equation}
\end{definition}

\vspace{-.2in}

If $\phi(\cdot,\cdot)$ is u\ped\ with respect to the whole state $x$ then we will simply say that ``$\phi$ is u\ped''. This notation will allow us to establish some results for nonlinear systems with state $x$ by imposing, on a certain function, the condition of  u\ped\ w.r.t. only $part$ of the state. 

For the sake of comparison we recall below the  definition proposed in \cite{TACDELTAPE} and an equivalent characterization of u\ped. The former is stated as a property of a pair of functions $(\phi,F)$ where $F$ is the vector field in \rref{one} and the {\em matrix} function $\phi:\mR\times\mR^n\to\mR^{p\times q}$ is such that  $\phi(\cdot\,,x(\cdot\,,\t\circ,\x\circ))$ is locally integrable for each solution $x(\cdot\,,\t\circ,\x\circ)$ of \rref{one}. The latter may be viewed as a trajectory-independent formulation of the former.

\vspace{-.1in}

\begin{definition}\cite{TACDELTAPE}
  \label{def:olddpe}
The pair $(\phi,\,F)$ is called uniformly $\delta$-persistently exciting (u$\delta$-PE) with respect to $\x1$ [\,along the trajectories of \rref{one}\,] if, for each $r$ and $\delta >0$ there exist constants $T(r,\delta) > 0$ and $\mu(r,\delta) > 0$,  s.t. $\forall (\t\circ,\x\circ)\in \mR\times {\cal B}(r)$, all corresponding solutions satisfy      
\begin{eqnarray*}
&\dty \left\{ \min_{\tau \in [t,\, t+T]} \norm{\x1(\tau)} \geq \delta \right\} \
\Rightarrow &\nonumber \\
&\dty \left\{
 \int_t^{t+T}
\phi(\tau,x(\tau,\t\circ,\x\circ))   \phi(\tau,x(\tau,\t\circ,\x\circ))^\top
d\tau \, \geq \mu  I \right\}\,&
\end{eqnarray*}
for all  $t \geq \t\circ$. 
\end{definition}
We emphasize  that Definition \ref{def:olddpe} is cited here only for the sake of comparison. Throughout this paper, when we say that a function is u\ped\ we mean in the sense of Definition \ref{def:newdpe}. 

Even though in spirit the properties in both definitions are the same, as pointed out before, they are technically different. This is illustrated by the following example.
\begin{example}\label{eg:dpedefs}
Consider the system 
\begin{equation}
\dot x := 
\begin{bmatrix}
\dot x_1\\
\dot x_2 
\end{bmatrix} = 
\begin{bmatrix}
-x_2\\
 x_1
\end{bmatrix}
=: f(t,x) 
\end{equation}
whose solutions with initial conditions $t_\circ=0$ and $x_\circ=[1\, 0]^\top$ take the form $x_{1}(t) = \cos t$, $x_{2}(t)=\sin(t)$.
Consider also the function 
\[
\psi(t,x):= [\sin (t) \, -\cos (t) ]
\begin{bmatrix}
x_1\\ x_2\,
\end{bmatrix}\,.
\]
Clearly, $\psi(\cdot,x)$ is locally integrable for each $x\in\mR^2$. One can see also that this function satisfies Definition \ref{def:newdpe} since $\norm{\psi(t,x)} = (x_1^2\sin t^2 -2x_1x_2\sin (t)\cos (t) + x_2^2\cos (t)^2)^{1/2}$ but it does not satisfy  Definition \ref{def:olddpe} with the initial conditions above, since $\norm{\psi(t,x(t))} \equiv 0$. On the other hand, the function 
$ \psi(t,x):=  x_2 \, $
satisfies the trajectories-dependent property of Definition \ref{def:olddpe} but it does not satisfy Definition \ref{def:newdpe} for any $x$ such that $x_1\neq 0$ (in particular for $x = [1\, 0]^\top$). However, it is u\ped\ in the sense of Def.  \ref{def:newdpe} with respect to $x_2$.
\end{example}

\subsection{Characterizations of u\ped}

We present now some useful  properties which are equivalent to Definition \ref{def:newdpe}. Our first ``characterization'' of u\ped\ is actually a relaxed property. It states that when dealing with the particular (but fairly wide) class of uniformly continuous functions, it is sufficient to verify the integral in \rref{eq:newdpe} only for each fixed $x$ such that $x_1\neq 0$  (i.e., for ``large'' states).
\begin{lemma}
\label{lem:equiv:defs} If $x\mapsto \phi(t,x)$ is continuous uniformly in $t$ then $\phi(\cdot,\cdot)$ is u\ped\ with respect to $x_1$ if and only if 
\begin{statement}
\label{stmt:B}  for each $x\in{\cal D}_1$ there exist $T>0$ and $\mu>0$ such that, for all $t\in\mR$, 
\begin{equation}\label{eq:B}
    \int_{t}^{t+T} |\phi(\tau,x)| d \tau \geq \mu \ .
\end{equation}
\pullup{12mm}
\end{statement}
\end{lemma}
The following Lemma helps us to see that both Definitions \ref{def:olddpe} and \ref{def:newdpe} state in words that ``a function $\phi(t,x)$ is u\ped\ with respect to $x_1$ if $t\mapsto \phi(t,x)$ is PE in the usual sense\footnote{That is, as defined for functions which depend only on time: that the function $A:\mRp\to\mR^{m\times n}$, $m\leq n$ is PE if  there exist $T>0$ and $\mu>0$ such that for all unitary vectors $z\in\mR^n$ we have that $\int_t^{t+T} z^\top A(s)^\top A(s) z\, ds \, \geq\, \mu$. } whenever the states $x_1$ (or similarly, the trajectories $x_1(t)$) are large''. This is important since it is the central idea to keep in mind when establishing sufficiency results based on the u\ped\ property. This idea also establishes a relation with the original but also technically different definition given in \cite{uped}. 
\begin{lemma}\label{TA:equiv2:B} The function $\phi(\cdot,\cdot)$ is u\ped\ w.r.t. $x_1$ if and only if 
\begin{statement}
\label{stmt:A} for each $\delta>0$ and $\Delta \geq \delta$ there exist $T>0$ and $\mu>0$ such that, for all $t\in\mR$, 
\begin{equation}\label{eq:E}
\dty \norm{x_1} \in [\delta,\,\Delta]\,, \ \norm{x_2} \in [0,\,\Delta]\qquad \Longrightarrow \qquad
\dty \int_t^{t+T} \norm{\phi(\tau,x)}d\tau \geq \mu\,. 
\end{equation}
\end{statement}\pullup{9mm}
\end{lemma}
The last characterization is useful as a technical tool in the proof of convergence results as we will see in Section \ref{sec:suf}.
\begin{lemma}\label{lem:equiv:PE}
The function $\phi(\cdot,\cdot)$ is u\ped\ w.r.t. $x_1$ if and only if\\[-7mm]
\begin{statement}
\label{stmt:C} for each $\Delta>0$ there exist $\gamma_\Delta\in\cK$ and $\theta_\Delta :\mR_{>0}\to\mR_{>0}$ continuous strictly decreasing such that, for all $t\in\mR$,
\begin{equation}\label{eq:equiv:PE}
\dty\left\{\, \norm{x_1},\ \norm{x_2} \in [0,\,\Delta]\backslash\{x_1=0\}\,\right\}  \qquad \Longrightarrow \qquad \dty \int_{t}^{t+\theta_\Delta(\norm{x_1})} \norm{\phi(\tau,x)} d\tau \geq \gamma_\Delta(\norm{x_1})\,.
\end{equation}
\end{statement}
\end{lemma}



\begin{remark}
We may summarize the above characterizations as follows.
\begin{itemize}
\item  The following are equivalent: $\phi$ is u\ped\ with respect to $x_1$, statement {\bf \ref{stmt:A}}, statement {\bf\ref{stmt:C}}. Also, each of these implies Statement {\bf\ref{stmt:B}}.
\item For  uniformly continuous functions, $\phi$ , it is sufficient to check Statement {\bf\ref{stmt:B}} which on occasions may be easier to verify. Consequently, {\bf\ref{stmt:A}}, {\bf\ref{stmt:C}} will also hold.
\end{itemize}
\end{remark}

\subsection{Properties and facts of u\ped\ functions}

We present below some important properties which hold for functions which are u\ped\ in the sense of Definition \ref{def:newdpe}. Similar properties are well known for only-time dependent functions (cf. e.g. \cite{NARANA}) and some of which were proved to hold  as well in \cite{TACDELTAPE}, for pairs $(F,\phi)$ which satisfy Definition \ref{def:olddpe}.

Let $\phi(t,x)\in\mR^m$ be u\ped\ with respect to $x_1$ in the sense of Definition \ref{def:newdpe} and let  $\alpha_i : \mRp\to\mRp$, $i=1,\dots, 2$  be continuous non decreasing functions. Assume that for all $x\in\mR^n$ and almost all $t\in\mR$,
\begin{eqnarray}
\label{eq:bndphi}
\norm{\phi(t,x)} & \leq & \alpha_1(\norm{x_1}) \\
\label{eq:bnddphi}
\norm{\jac{\phi(t,x)}{t}} & \leq & \alpha_2(\norm{x_1})\,.
\end{eqnarray}

We start by pointing out the following useful observation.
\begin{fact}[Power of a u\ped\ function]\label{prop:power}
If the {\em scalar} function $\phi(t,x)$ is u\ped\ with respect to $x_1$, with parameters $\delta$, $T$ and $\mu>0$ then, for any $p>1$ the function $\Psi(t,x):= \phi(t,x)^p$ is u\ped\ with respect to $x_1$, with parameters $T_p =T $ and $\mu_p := \frac{\mu^p}{T^{p/q}}$, $\frac{1}{p} + \frac{1}{q} =1$.
\end{fact}
\begin{claimproof}
This property follows straightforward using Cauchy-Schwartz inequality  to see that for each $p>1$ there exists $q$ satisfying $\frac{1}{p} + \frac{1}{q} =1$ and
\[
\int_{t}^{t+T} \norm{\phi(t,x)} dt\leq \left(\int_{t}^{t+T} \norm{\phi(t,x)}^pdt\right)^{1/p} T^{1/q}
\]
hence, for each $x\in {\cal D}_1$ and all $t\in\mR$ we have that 
\[
 \norm{z-x} \leq \delta \ \Rightarrow \  \int_{t}^{t+T} \norm{\phi(t,x)}^p dt \geq \frac{\mu^p}{T^{p/q}}\,.
\]
\end{claimproof}

It is also useful to remark that for functions linear in the state, u\ped\ is equivalent to the usual PE property (when restricting $t$ to the set of non-negative reals).
\begin{fact}\label{dpe:linearcase}
Consider the function $\phi(t,x):=\Phi(t)^\top x$ where $\Phi:\mR\to\mR^{n\times m}$ is locally integrable. Then, $\phi(t,x)$ is u\ped\ with respect to $x$ if and only if there exist $T$ and  $\mu>0$ such that  
\begin{equation}
\label{usualpe}
\int_{t}^{t+T} \Phi(\tau) \Phi(\tau)^\top d\tau \geq \mu\qquad \forall \,t\in\mR \,.
\end{equation}
\end{fact}
\begin{proof}
Sufficiency follows from Lemma \ref{TA:equiv2:B}. Necessity follows from the observation that the implication \rref{eq:E} holds in particular for $\delta=\Delta=1$ and that in this case, we have that for all $t\in\mR$ and for all unitary vectors $x\in\mR^n$,
\begin{equation}
\label{1000}
\int_{t}^{t+T} x^\top\Phi(\tau) \Phi(\tau)^\top x\, d\tau \geq \mu \,.
\end{equation}
This is well known\footnote{Strictly speaking this fact is well known for $t\in\mRp$ but clearly, this implies that the same property follows for $t\in\mR$.} to be equivalent to \rref{usualpe} (see e.g. \cite{MORNAR2}). The sufficiency follows observing that under \rref{usualpe}, the integral in \rref{1000} is lower-bounded by $\mu\norm{x}^2$ which in particular, for each pair $\Delta\geq \delta>0$ and for all $x$ such that $\Delta\geq \norm{x}\geq \delta$, is equal to  $\mu\delta^2=:\mu_\delta>0$.
\end{proof}

\begin{property}[Filtered u\ped\ function]\label{prop:fildpe}
Consider the differential equation
\begin{equation}
\label{eq:dphif}
\dot \Phi_f = -f_\phi(t,\Phi_f) \Phi_f + \phi(t,z)
\end{equation}
where  $\phi(t,z)$ is u\ped\ with respect to $z_1$ and assume that: \\
1) $f_\phi$ is such that 
\begin{equation}
\label{eq:bndfphi}
\norm{f_\phi(t,\Phi_f)} \leq \alpha_3(\norm{\Phi_f}),\quad \alpha_3\in {\cal N}\,.
\end{equation}
2) There exist $\alpha_4$, $\alpha_5\in {\cal N}$ such that with $\tilde \phi_f(\cdot\,,t_{f_\circ},\phi_{f_\circ},z)$ denoting the solution of \rref{eq:dphif}, we have that
\begin{equation}\label{eq:bndphift}
\norm{\tilde\phi_f(t,t_{f_\circ},\phi_{f_\circ},z)} \leq \alpha_4(\norm{\phi_{f_\circ}})\alpha_5(\norm{z_1})\quad \forall\,t\geq t_{f_\circ}\,.
\end{equation}
Then, defining $x_2:=\col[t_{f_\circ},\phi_{f_\circ},z_2]$ and $x_1:=z_1$, the function $\phi_f(t,x):=\tilde\phi_f(t,t_{f_\circ},\phi_{f_\circ},z)$ is u\ped\ with respect to $x_1$. 
\end{property}
\begin{claimproof}
We prove this property by showing that Statement {\bf\ref{stmt:A}} holds. Defining $\omega(t,x):= -\phi_f(t,x)^\top\phi(t,z)$  we have that 
\[
\dot \omega(t,x) = -\phi_f(t,x)^\top \jac{\phi}{t} + \big[ f_\phi(t,\phi_f)\phi_f - \phi(t,z)\big]^\top \phi(t,z)
\]
hence, dropping the arguments for notational simplicity, we get from \rref{eq:bnddphi} and \rref{eq:bndfphi} that 
\begin{equation}
\label{eq:dtomega}
\dot{ \omega} \leq\norm{\phi_f}\alpha_2(\norm{x_1}) + \alpha_3(\norm{\phi_f})\norm{\phi_f}\norm{\phi}- \phi^\top\phi\,.
\end{equation}
Then, using \rref{eq:bndphi} and \rref{eq:bndphift} we obtain that 
\[
\dot{ \omega} \leq\norm{\phi_f}\alpha_2(\norm{x_1}) + \alpha_3\big( \alpha_4(\norm{\phi_{f_\circ}})\alpha_5(\norm{x_1})\big)\norm{\phi_f}\alpha_1(\norm{x_1}) - \phi^\top\phi
\]
therefore, for each $\Delta > 0$ we have that for all $x_1$ and $\phi_{f_\circ}$ such that $\max\left\{\norm{\phi_{f_\circ}},\,\norm{x_1} \right\}  \leq \Delta$, the constant
\[
c :=  \alpha_2(\Delta) + \alpha_3\big(  \alpha_4(\Delta)\alpha_5(\Delta)\big)\alpha_1(\Delta) > 0 
\]
Inverting the sign and integrating on both sides of \rref{eq:dtomega} from $t$ to $t+T_f$  for some $T_f > T >0$ we obtain that
\begin{equation}
\label{eq:dtomega2}
\omega(t,x)-  \omega(t+T,x) \geq -c\int_t^{t+T_f} \norm{\phi_f(\tau,x)} d\tau + \int_t^{t+T_f} \norm{\phi(\tau,z)}^2d\tau\,.
\end{equation}
Since $\phi$ is u\ped\ with respect to $z_1=x_1$ with parameters $\mu$ and $T$ we can  invoke Fact \ref{prop:power} to lower bound the last term on the right hand of \rref{eq:dtomega} and, \rref{eq:bndphi}, \rref{eq:bndphift} to bound the terms on the left hand side. We obtain that for each $\Delta>\delta>0$, each $x_1$ and $\phi_{f_\circ}$ such that $\delta \leq\norm{x_1}   \leq \Delta$ and $ \norm{\phi_{f_\circ}}  \leq \Delta$,
\[
\int_t^{t+T_f} \norm{\phi_f(\tau,x)} d\tau \geq -\frac{2\alpha_4(\norm{\phi_{f_\circ}})\alpha_5(\norm{x_1})\alpha_1(\norm{x_1})}{c} + (k+1)\frac{\mu^2}{T}\,,
\]
where we have defined $k:= \lfloor T_f/T \rfloor - 1$ and $\lfloor r \rfloor$ as the largest  integer smaller than $r$. Thus, choosing $T_f(\Delta)$ sufficiently large so that 
\[
k \geq \frac{2\alpha_4(\Delta)\alpha_5(\Delta)\alpha_1(\Delta)T}{c(\Delta)\mu^2}
\]
we obtain that $\phi_f$ is u\ped\ with respect to  $x_1=z_1$ .
\end{claimproof}

\section{u\ped\ is necessary and sufficient for UGAS}

We present in this section our main results. We will show that for a fairly general class of nonlinear time-varying systems the property of u\ped\ is necessary and sufficient for uniform attractivity of the origin. In other words, we establish that for UGS systems u\ped\ is necessary and sufficient for UGAS.

\subsection{Necessity}
\label{sec:nec}

The following result, contained in \cite[Theorem 6.2]{ARST78}, gives
conditions under which u\ped\ of the right-hand side of a differential
equation is necessary for uniform asymptotic stability.
The technical conditions that we use permit a relatively straightforward
proof, based on Gronwall's lemma, without recourse to the notion of limiting equations which are used in \cite[Theorem 6.2]{ARST78}.

\begin{theorem}[UGAS $\Rightarrow$ u\ped]\label{thm:arstein}
Assume that $F(\cdot,\cdot)$ in \rref{one} is  Lipschitz in $x$ uniformly in $t$. If 
\rref{one} is UGAS, then $F(\cdot,\cdot)$ is u\ped\ with respect to $x\in \mR^n$.
\end{theorem}

\begin{proof}
The Lipschitz assumption on  $F(\cdot,\cdot)$ implies that for each $M>0$ such that $\max\{ \norm{y},\norm{z} \} \leq M$ there exists $L>0$ such that 
\begin{equation}
\label{eq:flip}
 \norm{F(t,y)-F(t,z)} \leq L\norm{y-z} \,
\end{equation}
for all $t\in\mR$.\remfootnote{This is true because $F(\cdot,\cdot)$ is uniformly locally Lipschitz for all $x$, in particular, it is locally Lipschitz on each {\em compact} ${\cal H}(0,M)$ which means that it is uniformly Lipschitz on ${\cal H}(0,M)$.} From the UGAS assumption on \rref{one} it follows that 
$\exists \beta\in\cKL$ s.t. $\forall (t_\circ,x_\circ)\in\mR\times\mR^n$
\begin{equation}
\label{eq:uasbeta}
\norm{x(t,t_\circ,x_\circ)} \leq \beta(\norm{x_\circ},t-t_\circ),\quad \forall \, t\geq t_\circ\,.
\end{equation}
Notice that without loss of generality, we may assume that $\beta(s,0) \geq s$ for all $s\geq 0$. 

For the purposes of establishing a contradiction, assume that the statement of the theorem does not hold. More precisely, and using Lemma \ref{lem:equiv:defs},
assume there exists $x^{*} \neq 0$ such that
 for each $T$ and $\mu>0$, there exist $t_*\in\mR$ such that
\begin{equation}
\label{eq:fnotdpe}
\int_{t_*}^{t_*+T} \norm{F(t,x_*)}dt < \mu\,.
\end{equation}
Pick $T$ such that $\beta(|x_*|,T) \leq \dty\frac{|x_*|}{2}$. 
Let the Lipschitz assumption generate the constant $L>0$ for 
$M:=\beta(|x_*|,0)$ and let $\mu:=\dty\frac{|x_*| e^{-LT}}{2}$. Let
these values generate $t^{*}$, and consider
the solution of \rref{one} starting at $(t^{*},x^{*})$.
By definition of solution,
\begin{equation}
\label{eq:solf}
x(t,t_*,x_*) = x_* + \int_{t_*}^t F(\tau,x(\tau))d\tau\,\quad \forall\, t\geq t_*\,,
\end{equation}
which also satisfies\remfootnote{ 
\begin{eqnarray}
\label{eq:solf2}
x(t,t_*,x_*) - x_*  & \leq &  \int_{t_*}^t F(\tau,x(\tau))d\tau\, +  \int_{t_*}^t F(\tau,x_*)d\tau -  \int_{t_*}^t F(\tau,x_*)d\tau\,\\
\label{eq:solf3}
\norm{x(t,t_*,x_*) - x_*}  & \leq &  \int_{t_*}^t \norm{ F(\tau,x(\tau)) - F(\tau,x_*)}d\tau\, + \int_{t_*}^t F(\tau,x_*)d\tau\,
\end{eqnarray}}
\begin{equation}
\label{eq:solf4}
\norm{x(t,t_*,x_*) - x_*}   \leq  \int_{t_*}^t L \norm{ x(\tau) - x_*}d\tau\, + \int_{t_*}^t \norm{F(\tau,x_*)}d\tau\,.
\end{equation}
Now, setting $t=t_*+T$ we have from \rref{eq:fnotdpe} that 
\begin{equation}
\label{eq:solf5}
\norm{x(t_*+T,t_*,x_*) - x_*}   <  \int_{t_*}^{t_*+T} L \norm{x(\tau) - x_*}d\tau\, + \mu\,
\end{equation}
so using the Gronwall-Bellman inequality we  obtain that 
\begin{equation}
\label{eq:solf6}
\norm{x(t_*+T,t_*,x_*) - x_*}   <   \mu e^{LT} = \frac{|x_*|}{2}\,.
\end{equation}
On the other hand, $\norm{x(t_*+T,t_*,x_*)}\leq \beta(\norm{x_*},T) \leq  
\dty\frac{|x_*|}{2}$ hence\remfootnote{So $-\norm{x(t_*+T,t_*,x_*)}\geq -\frac{\delta_*}{2}$ } 
\begin{eqnarray*}
\label{eq:solf7}
\norm{x(t_*+T,t_*,x_*) - x_*} &  \geq& \norm{x_*} - \norm{x(t_*+T,t_*,x_*)}  \\
\label{eq:solf8}
& \geq & \frac{|x_*|}{2}\,
\end{eqnarray*}
which  contradicts \rref{eq:solf6}. 
\end{proof}

\subsection{Sufficiency}

Our main results on sufficiency of u\ped\ for UGAS derive
from the sufficient conditions for UGAS that we have recently
established in \cite{MATCDC02}.  In that work we have generalized
Matrosov's theorem (see, for example, \cite[Theorem 5.5]{ROUMAW}).
In our generalization, we assume UGS, which is standard in Matrosov's
theorem, but then we use
an arbitrary finite number of auxiliary functions
to establish uniform convergence.  (Matrosov's theorem
relies on one auxiliary function.) Accordingly, the following
sufficient conditions are expressed in terms of a finite number
of auxiliary functions (Assumption \ref{assume:2})
having a certain nested property
(Assumption \ref{assume:3}) and
the property that when the bounds on the derivatives of these
auxiliary functions are all zero, this implies that part of
the state and a certain function are zero (Assumption
\ref{assume:4}) and, moreover, that the function is u\ped\ with
respect to the rest of the state (Assumption \ref{assume:5}).
The last technical assumption (Assumption \ref{assume:5b}) bounds
the derivative of the rest of the state in terms of the part of
the state and the function that are zero in Assumption \ref{assume:4}.
\label{sec:suf}
\begin{theorem}
\label{theorem:2}
Under Assumptions \ref{assume:1}-\ref{assume:5b} below, the origin of (\ref{one}) is UGAS.
\end{theorem}
\vspace{-.2in}
\begin{assumption}
\label{assume:1}
The origin is UGS.
\end{assumption}
\vspace{-.2in}
\begin{assumption}
\label{assume:2}
There exist integers $j$, $m > 0$ and for each $\Delta >0$  there exist
\begin{itemize}
\item a number $\mu >0$
\item locally Lipschitz continuous functions\\ $V_i:\mR\times\mR^n\to \mR$, $i\in \{1,\ldots,j \}$ 
\item a (continuous) function $\phi:\mR\times\mR^n\to \mR^{m}$,
\item continuous functions $Y_i:\mR^n\times\mR^m\to \mR$,\\ $i\in \{1,\ldots,j \}$
\end{itemize}
such that,
for almost all $(t,x) \in \mR\times \BD$, 
\begin{eqnarray}
\label{bndonViPhii}
&\dty \max\left\{ \norm{V_i(t,x)},\ \norm{\phi(t,x)} \right\}  \, \leq\, \mu, & \\
\label{bndondotVi}
&\dty \dot V_i(t,x) \,\leq\, Y_i(x,\phi(t,x)) \,.&  \ 
\end{eqnarray}
\end{assumption}
\begin{assumption}
\label{assume:3} 
For each $k \in \left\{1,\cdots,j\right\}$ we have that\footnote{For the case $k=1$ one should read Assumption \ref{assume:3} as $Y_1(x,\phi(t,x))\leq 0$ for all $(z,\psi) \in \BD \times \Bmu$. }
\begin{description}
\item{(A):}\ 
$\left\{ \quad  (z,\psi) \in \BD \times \Bmu  \ , \ Y_i(z,\psi) = 0 \right. $ $ \left. \forall i\in \{ 1,\ldots, k-1\}  \quad \right\}$
\end{description}
implies 
\begin{description}
\item{(B):}
$\quad \{ \ Y_k(z,\psi) \leq 0 \ \}$ .
\end{description}
\end{assumption}
\begin{assumption}
\label{assume:4} 
We have that
\begin{description}
\item{(A):}
$\quad \left\{ \ (z,\psi) \in \BD \times \Bmu \ , \quad
\ Y_i(z,\psi) = 0  \right. $ $ \left. \forall i\in \{ 1,\ldots, j\} \ \right\}$
\end{description}
implies
\begin{description}
\item{(B):}
$\quad \{ \ z_{1} = 0, \psi_{1} = 0 \ \}$ .
\end{description}
\end{assumption}
\begin{remark}
In Assumption \ref{assume:4},
there is no requirement that the size of $z_{1}$ matches the
size of $\psi_{1}$.
\end{remark}
\begin{assumption}
\label{assume:5}
The first component of $\phi(t,x)$ i.e., $\phi_{1}$, is independent of $x_{1}$, locally Lipschitz in $x_{2}$
uniformly in $t$, $u \delta$-PE w.r.t. $x_{2}$ and zero at the origin.
\end{assumption}
\begin{assumption}
\label{assume:5b} For all $(t,x) \in \Re \times \BD$, we have
$|F_{2}(t,x)| \leq \rho_{\Delta}(x_{1},\phi_{1}(t,x))$ where 
$\rho_{\Delta}$ is continuous
and zero at zero.
\end{assumption}
\begin{proof}[of Theorem \ref{theorem:2}]
This result follows from the main result in \cite{MATCDC02} by using the additional function
\begin{equation}
   V_{j+1}(t,x):= - \int_{t}^{\infty} e^{t-\tau} |\phi_{1}(\tau,x)| d \tau\,.
\end{equation}
We claim that this function satisfies for almost all $(t,x) \in \Re \times \BD$,
\begin{equation}\label{lego}
\shiftleft{3mm}   \begin{array}{ccl}
\dot{V}_{j+1}(t,x)     &\!\! \leq &\!\!   |\phi_{1}(t,x)| + Y_{a}(x_{2}) 
            + K(\Delta) \rho(x_{1},\phi_{1}(t,x))
   \end{array}
\end{equation}
where $$Y_{a}(x_{2}):= \max\{-\norm{x_2}\,,\,-\mbox{exp}(-\theta_\Delta(\norm{x_2}))\gamma_\Delta(\norm{x_2})\}\,$$
and $K(\Delta)$ represents a Lipschitz constant for
$\phi_{1}$ on the set $\Re \times \BD$.
 Then, the result follows defining $$Y_{j+1}(x,\phi(t,x)):= |\phi_{1}(t,x)| + Y_{a}(x_{2}) + K(\Delta) \rho(x_{1},\phi_{1}(t,x))\,.$$
 
We now prove that \rref{lego} holds. Let  the u\ped\ property of $\phi_1$ and Lemma \ref{lem:equiv:PE} generate the functions $\gamma_\Delta$ and $\theta_\Delta$. It is direct to see that for any $x \in {\cal B}(\Delta)\backslash\, \{x_2=0\}$ we have that 
\begin{eqnarray}
V_{j+1}(t,x) & \leq & -\mbox{exp}(-\theta_\Delta(\norm{x_2}))\int_{t}^{t+\theta_\Delta(\norm{x_2})} \norm{\phi_1(t,x)}d\tau \nn \\
\label{bndonVjp1}
& \leq & -\mbox{exp}(-\theta_\Delta(\norm{x_2}))\cdot\gamma_\Delta(\norm{x_2}) \, \leq\, 0\,.
\end{eqnarray}
Secondly, the partial derivatives of $V_{j+1}$ are given by
\begin{eqnarray*}
\jac{V_{j+1}}{t} &=&  \norm{\phi_1(t,x)} -\int_t^\infty \jac{}{t}\mbox{e}^{(t-\tau)}\norm{\phi_1(\tau,x)}d\tau\\
\jac{V_{j+1}}{x_2} &=& -\int_t^{\infty} \mbox{e}^{(t-\tau)}\jac{\norm{\phi_1(t,x)}}{x_2}d\tau\qquad \mbox{a.e.}
\end{eqnarray*}
The Lipschitz assumption on $\phi_1$ implies that for each $\Delta$ there exists $K(\Delta)$ such that, 
$$ \norm{\jac{\phi_1(t,x)}{x_2}} \ \leq \ K(\Delta) \qquad \mbox{a.e.}$$ 
Hence, we obtain that
\begin{eqnarray*}
\jac{V_{j+1}}{t} &  = & \norm{\phi_1(t,x)} + V_{j+1}(t,x)\\
\end{eqnarray*}
and for almost all $(t,x) \in \Re \times \BD$,
\begin{eqnarray*}
\jac{V_{j+1}}{x_2} & \leq & K(\Delta)\,. 
\end{eqnarray*}
It follows that for almost all $(t,x) \in \Re \times \BD$,
$$   \dot{V}_{j+1}(t,x)  \leq  |\phi_{1}(t,x)| + V_{j+1}(t,x) + K(\Delta) \rho(x_{1},\phi_{1}(t,x))\,.  $$
So \rref{lego} follows using the bound from \rref{bndonVjp1} in the inequality above. It is worth pointing out that $\theta_\Delta(\cdot)$ in \rref{bndonVjp1} is not defined at $x_2=0$ (see statement {\bf C} in Lemma \ref{lem:equiv:PE}). This motivates the definition of $Y_a(x_2)$ using the max of the two terms.
\end{proof}

The following corollary covers \cite[Theorem 6.3]{ARST78}.
\begin{corollary}
If the origin of (\ref{one}) is UGS and the following assumptions hold then, the origin of \rref{one} is also UGAS.
\end{corollary}
\begin{assumption}
\label{assume:7}
For each $\Delta>0$ there exist
\begin{itemize}
\item a number $\mu >0$
\item a locally Lipschitz continuous function\\ $V:\mR\times\mR^n\to \mR$,
\item a continuous function $Y:\mR^n\times\mR^n\to \mR$,\\ 
\end{itemize}
such that,
for almost all $(t,x) \in \mR\times \BD$, 
\begin{eqnarray}
\label{bndonViPhii:b}
&\dty \max\left\{ \norm{V(t,x)},\ \norm{F(t,x)} \right\}  \, \leq\, \mu, & \\
\label{bndondotVi:b}
&\dty \dot V(t,x) \,\leq\, Y(x,F(t,x))\,. &  \ 
\end{eqnarray}
\end{assumption}
\begin{assumption}
\label{assume:8} 
We have that
\begin{description}
\item{(A):}\ 
$\left\{ \quad  (z,\psi) \in \BD \times \Bmu  \right\}$
\end{description}
implies 
\begin{description}
\item{(B):}
$\quad \{ \ Y(z,\psi) \leq 0 \ \}$ .
\end{description}
\end{assumption}
\begin{assumption}
\label{assume:9} 
We have that
\begin{description}
\item{(A):}
$\quad \left\{ \ (z,\psi) \in \BD \times \Bmu \ , \quad
\ Y(z,\psi) = 0  \right\} $
\end{description}
implies
\begin{description}
\item{(B):}
$\quad \{ \psi = 0 \ \}$ .
\end{description}
\end{assumption}
\begin{assumption}
\label{assume:10}
$(t,x) \mapsto F(t,x)$ is locally Lipschitz in $x$ uniformly in $t$
and $u \delta$-PE with
respect to $x$.
\end{assumption}


\section{Applications in control systems}
\label{sec:appl}
\tonio{Maybe in this section we want to restrict time to the nonnegative reals? ;-)}

\subsection{Systems with matching nonlinearities}
\label{sec:MRAC}
We consider now a particular class of nonlinear  systems  \rref{one} with $x := \col[x_1, x_2]$, $x_1\in\mR^{n_1}$, $x_2\in\mR^{n_2}$  which may be viewed as a generalization of the classical strictly positive real system well studied in the context of linear systems. Hence, the structure of systems we consider in this section is important because, when inputs and outputs are considered\footnote{E.g. in applications of mechanical systems these may be external generalized torques and generalized velocities respectively.}, they naturally yield {\em passive} systems. Consequently, the type of results that we will present here may be used in the analysis of passivity-based (adaptive) control systems. This will be illustrated further below, with the application to adaptive control of robot manipulators. 

A typical example of systems with matching  nonlinearities  is what we may call {\em of the ``MRAC-type''} where MRAC stands for Model Reference Adaptive Control. These systems appear as closed loop equations in MRAC of {\em linear} plants (cf. \cite{MARTOM}). In the purely nonlinear context, they have also been for instance, in \cite{OrtegaFradkov,TACDELTAPE,JANIJACSP,THORSTUFFIJRNLC}. 

More precisely, we will consider systems of the form \rref{one} with
\begin{equation}\label{eq:spr}
F(t,x) := 
\begin{bmatrix}
A(t,x) + B(t,x)\\
C(t,x) + D(t,x)
\end{bmatrix}\,
\end{equation}
for which it is assumed that all functions are zero at $x=0$. Moreover, we will make the standing hypothesis that the system is UGS.

The reason why we decompose $F(t,x)$ in 4 terms is that we will impose different conditions on each of them. Roughly speaking, we will require that $\dot x_1 = A(t,x)$ is UGAS with respect to\footnote{We recall that a system is stable with respect to part of the state if, roughly speaking, the classical Lyapunov stability  properties hold for that part of the state. See \cite{VORIFAC02} for details.}  $x_1$, $A$ and $C$ vanish at $x_1 = 0$, $D$ and $B$ vanish at $x_2 = 0$ and moreover that the remainder of $B(t,x)$, i.e. when $x_1=0$, is u\ped.

The hypothesis on UGS  holds for a large class of systems, including  systems with the following interesting structure:
\begin{subeqnarray}\label{eq:error}
\dot{x}_{1} &  = & \tilde A(t,x_{1})x_1+G(t,x)x_{2}\label{eq:error1}\\
\dot{x}_{2} &  = & -P^{-1}G(t,x)^\top \left(  \frac{\partial W(t,x_{1})}{\partial x_{1}%
}\right)^\top ,\text{ \ \ \ \ }P=P^\top >0\label{eq:error2}%
\end{subeqnarray}
where $\x1\in \mR^{n_1}$, $\x2\in\mR^{n_2}$ and $W: \mR^{n_1}\times \mRp\to \mRp$ is a
${\cal C}^{1}$ positive definite radially unbounded function such that 
\begin{equation}
\label{WA}
\frac{\partial W(t,x_{1})}{\partial x_{1}} \tilde A(t,x_1)x_1 \leq 0\,.
\end{equation}
Indeed, it is sufficient to take $V(t,x):= W(t,x_1) + 0.5 x_2^\top P x_2$ to see that the system is UGS since $\dot V(t,x)\leq 0$. Notice that for this inequality to hold it is instrumental that the nonlinearities in the $x_2$-equation match with the second term in the $x_1$-equation. This motivates the title of the subsection.

It may be also apparent that the structure \rref{eq:error} is roughly, a  direct generalization of linear positive real systems and strictly positive real systems in the case when \rref{WA} holds with a bound of the form $-\alpha(\norm{x_1})$, with $\alpha\in \cK$. To better see this, let us  restrict $W(t,x_1)$ to be  quadratic then,  we can view \rref{eq:error} as two passive systems interconnected through the nonlinearity $G(t,x)$, i.e.,  the $x_1$ equation defines a passive map $ x_2\mapsto x_1$ (output feedback passive\footnote{See \cite{SEPBOOK} for a precise definition} if \rref{WA} holds with $-\alpha(\norm{x_1})$, with $\alpha\in \cK$) and the $x_2$ equation is an integrator (hence passive). See \cite{THOROBS} for a real-world example and \cite{MARTOM,THORSTUFFIJRNLC} for further discussions.

We address two cases: when $D(t,x_2)\equiv 0$ and when $ D(t,x_2) \not\equiv 0$. With reference to the observations above, if we regard $x_2$ as an input, we restrict $G(t,x)$ to depend only on $x_1$, and regard $\dot x_2$ as an output, these two cases actually correspond to those of relative degrees $1$ and $0$ respectively. This may be more clear if we restrict further our attention to linear time-varying systems and  define: $W(t,x_1):= \frac{1}{2}\norm{x_1}^2$, $z:=x_1$, $\dot z = A(t)z + B(t)u$ with $B(t)=G(t)$, $u:= x_2$, and output $y:=C(t)z+ \tilde D(t)u$  with $C(t):=P^{-1}B(t)$ and $\tilde D(t)u =: D(t,u)$. 

We also present some concrete examples in the following subsections: we first revisit some stabilization results for feedforward systems {\em a la} \cite{MAZPRA} and then, we see how our results apply to closed-loop identification (or adaptive control) of mechanical systems.

We now present the main results of this section. To that end, let us define
\begin{equation}
\label{def:Bo}
B_\circ(t,x_2) := B(t,x)\big|_{x_1=0}\,
\end{equation}
and notice that necessarily, $B_\circ(\cdot,0)\equiv 0$. 

We also introduce the following hypothesis which together with \rref{WA}, roughly speaking,  is related to the attractivity of the set $\{ x_1 = 0\}$ or in other words, to the inherent stability of $\dot x_1 = A(t,x)$ with respect to $x_1$ which, in particular implies that $x_1(t) \to 0$ as $t\to \infty$. Notice also that in the case that $A$ in \rref{eq:error} is linear time-independent and under \rref{WA}, the following assumption is equivalent to requiring that $\tilde A$ is stable.
\begin{assumption}
\label{assume:12} 
For the system defined by \rref{eq:spr} Assumption \ref{assume:7} holds and for the function $Y(\cdot,\cdot)$ in \rref{bndondotVi:b} we have that
\begin{description}
\item{(A):}
$\quad \left\{ \ (z,\psi) \in \BD \times \Bmu \ , \quad
\ Y(z,\psi) = 0  \right\} $
\end{description}
implies
\begin{description}
\item{(B):}
$\quad \{ x_{1} = 0 \ \}$ .
\end{description}
\end{assumption}

The next hypothesis imposes a regularity condition on $B(\cdot,\cdot)$ and that $A$ and $C$ vanish when $x_1=0$. Within the framework of linear adaptive control systems, we may say that the first part is a generalization of the global Lipschitz assumption on the regressor function. See \cite{TACDELTAPE} for further discussions.
\begin{assumption}\label{assume:13}
The functions $A$, $B$ and $C$ are locally Lipschitz in $x$ uniformly in $t$.
Moreover,  for each $\Delta\geq 0$ there exist $b_M > 0$ and continuous nondecreasing functions $\rho_i : \mRp \to \mRp$ such that $\rho_i(0)=0$ and for almost all $t\in \mR$ and  $x\in \mR^n$
\begin{equation}
\label{bM:bnd}
\max \,\left\{ \norm{B_\circ(t,x_2)},\,  \norm{\jac{B_\circ}{t}},\, \norm{\jac{B_\circ}{x_2}} \right\} \,\leq\, \rho_1(\norm{x_2})\,,
\end{equation}
and for all  $(t,x) \in \mR \times \mR^n$, $\norm{x_2}\leq \Delta$,
\begin{eqnarray}
\label{rho1:bnd}
&\dty \norm{B(t,x)-B_\circ(t,x_2)}  \leq \rho_2(\norm{x_1}) & \\[1mm] 
\label{rho2:bnd}
& \dty \max_{|x_2|\leq \Delta} \left\{ \norm{A(t,x)},\, \norm{C(t,x)} \right\}  \leq \rho_3(\norm{x_1}) &
\end{eqnarray}
\end{assumption}

The following theorem generalizes related results in the previously cited papers, including \cite{OrtegaFradkov,JANIJACSP,TACDELTAPE}. See the last reference for a detailed but non exhaustive literature review.
\begin{theorem}\label{thm:spr} 
Consider the system \rref{one}, \rref{eq:spr} under Assumptions \ref{assume:1}, \ref{assume:7}, \ref{assume:8}, \ref{assume:12} and \ref{assume:13}. Suppose also that 
\begin{assumption*}\label{ass:a5}
there exists a continuous non decreasing function $\rho_4:\mRp\to\mRp$ 
such that $\rho_4(0) = 0$, 
\begin{equation}
\label{bndonD}
 \norm{D(t,x)} \leq \rho_4(\norm{x_2})\,
\end{equation}
and
Statement {\bf\ref{stmt:C}} holds with $\theta_\Delta$ and $\gamma_\Delta$ such that for all $x_2\neq 0$
\begin{equation}\label{enoughpe}
\mbox{e}^{-\theta_\Delta(\norm{x_2})} \gamma_\Delta(\norm{x_2}) \geq 3\rho_1(\Delta)\rho_4(\norm{x_2})\,.
\end{equation}
\end{assumption*}
Then, the origin is UGAS.
\end{theorem} 
\begin{remark}
Notice that for the common case when $D\equiv 0$, Assumption \ref{ass:a5} reduces to requiring that $B_\circ(\cdot,\cdot)$ is u\ped\ with respect to  $x_2$. Also, in this case the necessity of  the latter follows directly  from Theorem \ref{thm:arstein} by observing that u\ped\ of $F(\cdot,\cdot)$ as defined in \rref{eq:spr} implies by virtue of \rref{rho2:bnd}, that $B_\circ(\cdot,\cdot)$ is u\ped\ with respect to $x_2$.
\end{remark}
\begin{proof}[of Theorem \ref{thm:spr}]
We  appeal to Theorem  \ref{theorem:2}. Assumption \ref{assume:1} is our standing hypothesis. To verify the rest of the assumptions of Theorem \ref{theorem:2} we introduce $V_{1}(t,x):=V(t,x)$ where $V(t,x)$ comes from Assumption \ref{assume:7} and the locally Lipschitz (due to Assumption \ref{assume:13}) function
\begin{equation}
\label{W}
V_{2}(t,x) = - x_1^\top B_\circ(t,x_2) - \int_t^\infty \mbox{e}^{(t-\tau)}\norm{B_\circ(\tau,x_2)}^2d\tau\,
\end{equation}
hence, $j=m=2$ in Assumption \ref{assume:2}. We also introduce 
\begin{equation}
\phi(t,x) :=
\begin{bmatrix}
x_2 \\ F(t,x)\,
\end{bmatrix}\,.
\end{equation}
With these definitions, it is clear that \rref{bndonViPhii} holds since all the functions are locally Lipschitz in $x$ uniformly in $t$ and $V_i(t,0)\equiv\phi(t,0)\equiv 0$. To see that \rref{bndondotVi} also holds, we observe first  that it is satisfied for $V_1(t,x)$ in view of \rref{bndondotVi:b} hence, we only need to check that there exists a continuous function $Y_2(\cdot,\cdot)$ satisfying the required conditions. To that end, we evaluate the total time derivative of $V_2(t,x)$ along the solutions of \rref{one}, \rref{eq:spr} to obtain\remfootnote{Here we are using that 
\[
\jac{}{t}\left[ - \int_t^\infty \mbox{e}^{(t-\tau)}\norm{B_\circ(\tau,x_2)}^2d\tau \right] = \norm{B_\circ(\tau,x_2)}^2 - \int_t^\infty \mbox{e}^{(t-\tau)}\norm{B_\circ(\tau,x_2)}^2d\tau - x_1^\top B_\circ(t,x_2) + x_1^\top B_\circ(t,x_2) 
\]
where the two middle terms equal to $V_2$. That is, the first term above corresponds to the last term in the middle line in \rref{dotW}.}
\begin{eqnarray}\nonumber
 \dot V_2(t,x) &=& -  x_1^\top\left(\jac{B_\circ}{t} +\jac{B_\circ}{x_2}[C(t,x) + D(t,x)]\right) + V_2(t,x) + x_1^\top B_\circ(t,x_2)\\
\nonumber
&&  - [A(t,x) + B(t,x)-B_\circ(t,x_2)]^\top B_\circ(t,x_2) \nonumber\\
&& -\left( \int_t^\infty \mbox{e}^{(t-\tau)} B_\circ(\tau,x_2)^\top \jac{ B_\circ(\tau,x_2)}{x_2}\right) [C(t,x) + D(t,x)]\,\quad \mbox{a.e}. 
 \label{dotW}
\end{eqnarray}

Notice also that by the u\ped\ condition on $B_\circ(\cdot,\cdot)$ and in view of Assumption \ref{assume:13} we have that 
\begin{equation}
\label{W:bnds}
V_2(t,x) \leq  \norm{x_1}\rho_1(\norm{x_2}) -  \mbox{e}^{-\theta_\Delta(\norm{x_2})}\gamma_\Delta(\norm{x_2})\,.
\end{equation}
Furthermore, using this and Assumption \ref{assume:13} again, we can over-bound several terms on the right hand side of \rref{dotW} as follows. Define $b_M:=\rho_1(\Delta)$ and $\bar\rho(r,s):= b_M[\, 3r + (r+3)\rho_3(r) + r\rho_4(s)  + \rho_2(r) \,]$ 
then, for almost all $(t,x)\in \mR\times {\cal B}(\Delta)$,\remfootnote{ 
\begin{eqnarray*}
\dot W(t,x) & \leq & \norm{x_1}[b_M + b_M(\,\rho_3(\norm{x_1})+\rho_4(\norm{x_2})\,)\,] +  \norm{x_1}\rho_1(\norm{x_2}) -  \mbox{e}^{-T_\Delta(\norm{x_2})}\theta_\Delta(\norm{x_2}) + \norm{x_1}b_M + \\
&& [\rho_3(\norm{x_1}) + \rho_2(\norm{x_1})]b_M + 2b_M\rho_3(\norm{x_1}) + 2b_M\rho_4(\norm{x_2})
\end{eqnarray*}}
\begin{equation}
\label{dotW:bnds}
\dot V_2(t,x) \leq \bar\rho(\norm{x_1},\norm{x_2})  -  \mbox{e}^{-\theta_\Delta(\norm{x_2})}\gamma_\Delta(\norm{x_2}) + 2b_M\rho_4(\norm{x_2})=: \widebar Y_2(x,\phi(t,x))\,.
\end{equation}
Define $Y_2(x,\phi(t,x)) := \max\{-\norm{x_2},\, \widebar Y_2(x,\phi(t,x)) \}$. In view of \rref{enoughpe} and  since $\bar\rho(0,s)\equiv 0$ we have that $Y_2(x,\phi(t,x)) \leq \max\{-\norm{x_2},\, - \frac{1}{3} \mbox{e}^{-\theta_\Delta(\norm{x_2})}\gamma_\Delta(\norm{x_2})\}\leq 0$ when $x_1 = 0$. Thus, Assumption \ref{assume:3} holds for $k=1$ due to Assumption \ref{assume:8} (we recall that here, $Y_1 = Y$) and for $k=2$, because $Y_2  (x,\phi(t,x)) \leq 0$ when $x_1 = 0$. 

Assumption \ref{assume:4} holds due to the following. Let $Y_1  (x,\phi(t,x)) = Y_2  (x,\phi(t,x)) =0$. Then,  by Assumption \ref{assume:12}   we have that $x_1 = 0$ while by definition, $Y_2(x,\phi(t,x)) =0$ implies that $x_2=0$.

Assumption \ref{assume:5} trivially holds and Assumption \ref{assume:5b} holds with $\rho_\Delta(x_1,x_2):= \rho_4(\norm{x_2}) + \rho_3(\norm{x_1})$.
\end{proof}

In support to the discussion at the beginning of the section it is worth remarking that Assumptions \ref{assume:1}, \ref{assume:7}, \ref{assume:8} and \ref{assume:12} hold under the following more restrictive but  commonly satisfied hypothesis, at least for a large class of {\em passive} systems. Below, we  present a concrete example concerning the adaptive tracking control of mechanical systems.
\begin{assumption}
\label{assume:11bis}
There exists a locally Lipschitz function $V:\Re \times \Re^{n} \rightarrow
\Re^{n}$, class-${\mathcal K}_{\infty}$ functions $\alpha_{1},\alpha_{2}$
and a continuous, positive definite function $\alpha_{3}$ such that
\begin{equation}\label{boundonVbis}
    \alpha_{1}(|x|) \leq V(t,x) \leq \alpha_{2}(|x|)
\end{equation}
and, almost everywhere,
\begin{equation}\label{boundondotVbis}
   \dot{V}(t,x) \leq - \alpha_{3}(x_{1}) \ .
\end{equation}
\end{assumption}
It may be more clear from this assumption that part of the conditions in Theorem \ref{thm:spr} are in the spirit of imposing that the system be asymptotically stable with respect to the $x_1$-part of the state. Then, for the convergence of the $x_2$-part of the trajectories we impose the u\ped\ exciteness condition. The following result which covers a class of systems similar to those covered by the main results in \cite[Appendix B.2]{MARTOM}, \cite{OrtegaFradkov,TACDELTAPE} is a direct corollary of Theorem \ref{thm:spr}.
\begin{proposition}\label{prop:ugas:MRAC:2}
The system \rref{one}, \rref{eq:spr} with $D\equiv 0$ and under Assumptions \ref{assume:13}, \ref{assume:11bis} is UGAS if and only if $B_\circ(\cdot,\cdot)$ is u\ped\ w.r.t.  $x_2$\,. 
\end{proposition}

For  clarity of exposition we present separately the following result, for systems  when the ``feedthrough'' term $D(t,x)$ is present.
\begin{proposition}\label{prop:spr:d}
Consider the system \rref{one}, \rref{eq:spr} under Assumptions \ref{assume:13}, \ref{assume:11bis} and \ref{ass:a5}. 
Then, the origin is UGAS.
\end{proposition}

We close the section with the counterparts of Propositions \ref{prop:ugas:MRAC:2} and \ref{prop:spr:d} which establish sufficient conditions for ULES. 

\begin{proposition}\label{prop:ules:1}
Consider the system \rref{one}, \rref{eq:spr} under Assumptions \ref{assume:13},  \ref{ass:a5},  and \ref{assume:11bis} and assume further that $\alpha_i(s):=\alpha_i\, s^2$, $i=1\ldots 3$, $\rho_j(s)=\rho_j\, s^{}$, $j=1\ldots 4$ for small $s$ and assume that the function $B_\circ(t,x_2)$ is u\ped\ in the sense that Statement {\bf\ref{stmt:C}} holds with functions $\gamma_\Delta(s)$ and $\theta_\Delta(s)$ such that $\mbox{e}^{-\theta_\Delta(s)}\gamma_\Delta(s) \geq \mu\,s^2$ for small $s$ and $\mu \geq 3b_M\rho_4$. Then, the origin is ULES. 
\end{proposition}
\begin{proposition}\label{prop:ules:2}
Consider the system  \rref{one}, \rref{eq:spr} with  $D(t,x)\equiv 0$ and let Assumptions \ref{assume:13},  \ref{ass:a5} and \ref{assume:11bis} hold with $\alpha_i(s):=\alpha_i\, s^2$, $i=1\ldots 3$, $\rho_j(s)=\rho_j\, s^{}$, $j=1\ldots 3$  for small $s$ and assume that the function $B_\circ(t,x_2)$ is u\ped. Then, the origin is ULES. 
\end{proposition}
\begin{proof}[of Propositions \ref{prop:ules:1} and \ref{prop:ules:2}]
We provide a combined proof for both propositions, based on standard Lyapunov theory.

Let $\Delta$ be generated by the u\ped\ assumption on $B_\circ(t,x_2)$. Let $0<R\leq\Delta$ be such that $\rho_i(s)=\rho_i\, s^{}$, $\alpha_i(s):=\alpha_i s^2$ and $\mbox{e}^{-\theta_\Delta(s)}\gamma_\Delta(s) \geq \mu\,s^2$ for all $s\leq R$. Let $r:=R\sqrt{\frac{\alpha_1}{\alpha_2}}$ then, from \rref{boundonVbis} and \rref{boundondotVbis} we obtain that $\norm{x_\circ}\leq r$ implies that $\norm{x(t)}\leq R$. In the sequel, we will restrict the initial conditions to $x_\circ\in {\cal B}(r)$.

Consider the Lyapunov function candidate ${\cal V}(t,x) := V(t,x) + \ep V_2(t,x)$ where $V_2(t,x)$ is defined in \rref{W} and $\ep$ is a small positive number to be chosen. Notice that $V_2(t,x)$ satisfies on $\mR\times {\cal B}(R)$, 
\begin{equation}
\label{W:bnds:ules}
-\ep\rho_1\norm{x_2}^2 - \ep \rho_1\norm{x_1}\norm{x_2}^{} \leq \ep V_2(t,x) \leq \ep \rho_1\norm{x_1}\norm{x_2}^{} - \ep \mu \norm{x_2}^2\,.
\end{equation}
 So we have that for sufficiently small $\ep$ there exist $\alpha_1'>0$ and $\alpha_2'>0$ such that for all $(t,x) \in \mR\times {\cal B}(R)$,
\begin{equation}
\label{bndsonVW}
\alpha_1'\norm{x}^2 \leq {\cal V}(t,x)  \leq \alpha_2'\norm{x}^2\,.
\end{equation} 
The total time derivative of ${\cal V}(t,x)$ on the points of existence and along the systems solutions yields, using \rref{boundondotVbis} and \rref{dotW:bnds},
\begin{equation}\label{fifty4}
\dot {\cal V}(t,x) \leq -(\alpha_3-\ep\nu_R)\norm{x_1}^2 + \ep\nu_R\norm{x_1}\norm{x_2}^{} - \ep\mu\norm{x_2}^2 + 2\ep b_M\rho_4\norm{x_2}^2\,\quad \mbox{a.e.}
\end{equation}
where $\nu_R>0$. Next, using the condition imposed on $\mu$ we obtain that
\begin{equation}\label{funfundfunfzig}
\dot {\cal V}(t,x) \leq -(\alpha_3-\ep\nu_R)\norm{x_1}^2 + \ep\nu_R\norm{x_1}\norm{x_2}- \ep b_M\rho_4\norm{x_2}^2\quad \mbox{a.e.}
\end{equation}
In the case that $D\equiv 0$ we have that $\rho_4 = 0$ and completing squares in \rref{fifty4} we obtain    that  for sufficiently small $\ep$ and sufficiently large $\alpha_3$ there exists $c>0$ such that $\dot {\cal V}(t,x) \leq -c\norm{x}^2$. Otherwise, the latter holds from \rref{funfundfunfzig} under similar arguments. ULES follows invoking standard Lyapunov theorems.

\end{proof}

\subsection{Dynamic bounded state feedback stabilization}
\def\umax{u_{\mbox{\scriptsize max}}}
We briefly revisit the problem of stabilizing a nonlinear time-varying system by a dynamic state feedback control law $u(t,\xi)$ under the constraint that, given $\umax>0$ the input $\norm{u(t,\xi(t))}\leq \umax$  for all $t\geq t_\circ$. In particular, we will revisit the generalized ``Jurdevic \& Quinn'' approach proposed in \cite[Section II.C, item 2)]{MAZPRA} for feedforward systems.

Consider the system 
\begin{equation}\label{mazpra:1}
\dot \xi = f(t,\xi) + g(t,\xi,u)u
\end{equation}
under the assumption that $f(t,\cdot)$ is locally Lipschitz uniformly in $t$ and $f(\cdot,x_2)$ is continuous,  $\dot \xi = f(t,\xi)$ is UGS and $g_\circ(t,\xi):=g(t,\xi,0)$ is such that for each $\Delta>0$ there exists $g_M>0$ such that 
\begin{equation}\label{bndsondg0}
\max_{\norm{\xi}\leq \Delta} \left\{\norm{g_\circ(t,\xi)},\,\norm{\jac{g_\circ(t,\xi)}{t}},\,\norm{\jac{g_\circ(t,\xi)}{\xi}}\right\} \,\leq\,g_M\,\qquad \mbox{a.e.}
\end{equation}
Under these conditions it is not difficult to show that there exist $\rho_f$ and $\rho_g\in {\cal N}$ such that $\norm{f(t,\xi)}\leq \rho_f(\norm{\xi})$ and the bound above holds with $\rho_g(\norm{\xi})$ for all $\xi$.

\begin{proposition}\label{prop:mazpra}
Assume that there exists a Lyapunov function  $V:\mRp\times\mR^n\to \mRp $ for  $\dot \xi = f(t,\xi)$ such that 
\[
\jac{V}{t} + \jac{V}{\xi}f(t,\xi) \leq 0
\]
and there exist two class $\cal K$ functions $\alpha_4$, $\alpha_5$ such that 
\begin{equation}\label{bndsondV}
\alpha_4(\norm{\xi}) \leq \norm{\jac{V}{\xi}} \leq \alpha_5(\norm{\xi})\,.
\end{equation}
The system \rref{mazpra:1} in closed loop with 
\begin{subeqnarray}\label{control:mazpra}
\dot z & = & -z - g(t,\xi,\Tanh(z))^\top \jac{V}{\xi}^\top \\
u & = & \Tanh(z), \qquad \Tanh(z):=\col[\tanh(z_1),\ldots,\tanh(z_n)]
\end{subeqnarray}
is UGAS  only if for any unitary vector $\zeta$, the function $g_\circ(t,\xi)^\top \zeta$ is u\ped. Furthermore, the origin is UGAS if this function is u\ped\ and Statement {\bf\ref{stmt:C}} holds with $\theta_\Delta$ and $\gamma_\Delta$ such that \rref{enoughpe} holds with $\rho_1(\Delta):=g_M$ and $\rho_4 = \rho_f$.
\end{proposition}
\begin{remark}
Certainly the result above holds for other (not necessarily smooth) saturation functions but $\tanh(\cdot)$ is particularly convenient for the Lyapunov analysis that follows.
\end{remark}
That is, with respect to the result in \cite[Section II.C, item 2)]{MAZPRA}  we contribute by addressing the problem for time-varying systems (see also \cite{MAZPRA99}) and relax 
 the ``observability-type'' condition \cite[Assumption A2]{MAZPRA} by imposing a u\ped\ ccondition. As a matter of fact, notice that we can actually deal with the case when $f(t,\xi)\equiv 0$ which is not included in the class of systems considered in \cite{MAZPRA}.

As an immediate corollary we also have the following. Consider the system 
\begin{equation}\label{mazpra:2}
\dot \xi = F(t,\xi,u)
\end{equation}
with  $F(t,\xi,\cdot)$ twice continuously differentiable and $\dot \xi = F(t,\xi,0)$ is UGS. Then, we can rewrite this system in the form \rref{mazpra:1} with 
\begin{eqnarray}
f(t,\xi ) & := & F(t,\xi,0)\\
g(t,\xi,u) & := & \int_0^1 \jac{F}{u}(t,\xi,su) ds\,.
\end{eqnarray}
Therefore, under the conditions above, the system \rref{mazpra:2} may be rendered UGAS by the bounded dynamic state feedback \rref{control:mazpra}.

\begin{proof}[of Proposition \ref{prop:mazpra}]
To see that the system is UGS we evaluate the total time derivative of the function\footnote{For other qualifying saturation functions we would use $\dty\sum_{i=1}^n\int_0^{z_i} \sat(s)ds$ in \rref{W:fwding}.}
\begin{equation}\label{W:fwding}
W(t,\xi,z) = V(t,\xi) + \sum_{i=1}^n \ln(\!\!\cosh(z_i))
\end{equation}
to obtain 
\begin{equation}
\dot W(t,\xi,z) \leq - z^\top \Tanh(z) \leq 0\,.
\end{equation}
To show that the system is uniformly attractive for such  initial conditions we simply apply Proposition \ref{prop:spr:d} with $x_1:=z$, $x_2 := \xi$, $A(t,x):= -z$, $\dty B(t,x):= -g(t,\xi,\Tanh(z))^\top\jac{V}{\xi}^\top$, $C(t,x):= g(t,\xi,z)z$ and $D(t,x):= f(t,\xi)$. For this, we observe that under \rref{bndsondV} we have that 
\begin{eqnarray*}
\norm{B_\circ(t,x_2)}^2 & = &  \jac{V}{x_2}g(t,x_2,0) g(t,x_2,0)^\top\jac{V}{x_2}^\top \\ 
&& \geq \norm{g(t,x_2,0)^\top \zeta } \alpha_4(\norm{x_2})^2 
\end{eqnarray*}
so it follows from Fact \ref{prop:power} and the monotonicity of $\alpha_3(\cdot)$  that $B_\circ(t,x_2)$ is u\ped\ if and only if $g(t,x_2,0)^\top \zeta = g(t,\xi,0)^\top \zeta$ is u\ped. The rest of the sufficient conditions of Proposition \ref{prop:spr:d} hold observing that $\sum_{i=1}^n \ln(\!\!\cosh(z_i))$ is positive definite and radially unbounded and letting $\alpha_3(x_1):= x_1^\top \Tanh(x_1)$.
\end{proof}
The next proposition is a direct corollary for systems without drift, i.e., with $f(t,\xi)\equiv 0$. This class of systems has been extensively studied from many viewpoints since the celebrated paper \cite{BRO}. While not explicitly mentioned except for recent references of the authors, a common point of many of the approaches based on time-varying controls for driftless systems are based on the implicit assumption that the control input ``excites'' the system in some sense. We may cite the following ``classical'' papers \cite{POM92,MURSAS,SAM1st,SAM3} and \cite{MORSAMEJC,JIANONHOL96} where the PE is implicitly used to prove convergence.  See also the more recent works \cite{uped,NONHOLEJC,SAMSIAM02}. The first two references present what we have called u\ped\ controllers. The last presents a general result on a characterization of controllability for driftless systems and covers as particular case, the problem of stabilization of nonholonomic systems.
\begin{proposition}[Driftless systems]
The system \rref{mazpra:1} with $f(t,\xi)\equiv 0$ and the controller \rref{control:mazpra}, is UGAS if and only if $g_\circ(t,\xi)^\top$ is u\ped. 
\end{proposition}
The proof follows straightforward along the lines of the proof of Proposition \ref{prop:mazpra} by taking $V(t,\xi)=\frac{1}{2}\norm{\xi}^2$.

\subsection{Closed loop identification of mechanical systems}

Lastly, we present a more particular example where the main results of this paper are useful. We will briefly address the problem of closed loop identification of mechanical systems, given by the Lagrangian equations:
\begin{equation}
\label{ELsystem}
D(q)\ddq + C(q,\dq)\dq + g(q) = u
\end{equation}
where all matrices and vector functions are smooth in their arguments.

We address the identification problem by considering the problem of tracking a desired (bounded and smooth) time-varying trajectory $\qd(t)$ under the assumption that positions and velocities $(q,\dq)$ are available for measurement and the system parameters are unknown. This problem has been studied and partially solved (i.e., guaranteeing the convergence of tracking errors) by many different approaches within the robot control community. See e.g. \cite{ORTSPO} for a tutorial on adaptive control and \cite{PBCELS} for a discussion from a  passivity viewpoint. See also the recent tutorial on identification of mechanical systems, \cite{GAUPOI}, including a broad list of references. 

We remark that the contents of this section were originally reported in \cite{ICRA03} and they extend the results of \cite{SLOTLIPE} where the authors proved non uniform parametric convergence under the {\em same} conditions that we will present here. Hence, among the numerous existing approaches to adaptive robot control we find that the best choice to highlight the utility of some of the results of this paper is the well known ``Slotine and Li'' algorithm \cite{SLOTLIA}.

To present this controller and a proof of UGAS of the origin (including the parametric errors) we recall that  the model \rref{ELsystem} is linear in the unknown parameters, i.e., there exists a smooth function $\Psi:\mR^n\times\mR^n\to\mR^{m\times n}$ such that
\begin{equation}
\label{ELsystem:psi}
D(q)\ddq + C(q,\dq)\dq + g(q) = \Psi(q,\dq,\ddq)^\top\theta\,.
\end{equation}
For particular examples of robot manipulators and the explicit expressions of $\Psi$ see for instance \cite{SPOVID}.

The well known Slotine and Li \cite{SLOTLIA} (passivity-based) adaptive control law is given by
\begin{equation}
\label{slotlicontrol}
u = \widehat D(q) \ddqr + \widehat C(q,\dq)\dqr + \hat g(q) - K_d s
\end{equation}
where $\hat {(\cdot)}$ denotes the estimate of the constant unknown lumped parameters, $\theta\in\mR^m$, contained in $(\cdot)$, $K_d$ is positive definite, $\dqr := \dqd - \lambda\tq$, with $\lambda > 0$, $\tilde q := q - q_*$ and $s := \dq - \dqr$. Again, from the property that the model is linear in the unknown parameters one can show that  there exists $\widetilde \Psi:\mR^n\times\mR^n\times\mR^n\times\mR^n\to\mR^{m\times n}$ such that, defining the parameter estimation error $\tilde \theta := \theta -\hat \theta$, the closed loop takes the form
\begin{equation}
\label{slotlicontrol:psi}
D(q)\dot s + C(q,\dq)s + K_d s = \widetilde \Psi(q,\dq,\dqr,\ddqr)^\top\tilde\theta\,.
\end{equation}
It is well known that the system above is uniformly globally asymptotically stable if $\widetilde \Psi^\top\tilde\theta\equiv 0$. For instance, in \cite{SPOORTKEL} a Lyapunov function with a negative definite bound on the total time derivative  is presented.

Consider as in many previous papers starting probably with \cite{SLOTLIPE}, that the estimated parameters vector $\hat\theta(t)$ is updated according with the so called {\em speed-gradient} adaptive law,
\begin{equation}
\label{adaptivelaw:a}
\dot{\hat \theta} = -\gamma \widetilde\Psi(q,\dq,\dqr,\ddqr) s\,, \qquad \gamma >0\,.
\end{equation}
Based on the theory for linear time-varying systems and classical results tailored for regressor functions which depend only on time, it was proved in \cite{SLOTLIPE} that the closed loop state trajectories tend to zero and, provided that $\widetilde \Psi(\qd(t),\dqd(t),\dqd(t),\ddqd(t))$ is persistently exciting, so do the parameter estimation errors.

Roughly speaking, this follows observing that $\widetilde \Psi(q(t),\dq(t),\dqr(t),\ddqr(t)) \, \longrightarrow\, \widetilde \Psi(\qd(t),\dqd(t),\dqd(t),\ddqd(t))$ as $t\to \infty$ and as it is well known now, (see e.g. \cite{NARANA}) a signal $\psi(t) = \psi^*(t) + \ep(t)$ where $\ep\to 0$ and $\psi^*$ is PE, remains PE.

The following proposition extends this result to UGAS of the origin of the closed loop system (that is, of $(s,\tilde q,\tilde \theta) = (0,0,0)$ ) with the {\em same} adaptive controller and under the {\em same} PE condition. Moreover, it also establishes that PE of $\widetilde \Psi(\qd(t),\dqd(t),\dqd(t),\ddqd(t))$ in the sense of \rref{usualpe} is {\em sufficient} and {\em necessary} for {\em uniform} parameter convergence. 

\begin{proposition}[Adaptive Slotine-and-Li controller \cite{ICRA03}]\label{prop:slotli}
Consider the system \rref{ELsystem} in closed-loop with \rref{slotlicontrol} and the  speed-gradient adaptive law \rref{adaptivelaw:a}. Assume that the reference trajectory $\qd(t)$ and its first two derivatives are bounded for all $t$. 
Then, the origin $(\tilde\theta,\tilde q, s)=(0,0,0)$ is UGAS if and only if the function 
$$\Phi(t):= \widetilde\Psi(\qd(t),\dqd(t),\dqd(t),\ddqd(t))$$
is persistently exciting in the sense of \rref{usualpe}. 
\end{proposition}
The proof (which we do not provide here) relies on Proposition \ref{prop:ugas:MRAC:2}. Roughly speaking, Assumption \rref{assume:13} holds from the smoothness and boundedness properties of the matrices $D(q)$, $C(q,\dq)$ and of the vector $g(q)$ (see any book on robot control, e.g. \cite{SPOVID,SCISIC}). Assumption \ref{assume:11bis} may be shown to hold using a similar Lyapunov function as in \cite{SPOORTKEL}. The most interesting is the u\ped\ condition. That this condition holds may be clear observing that in this case, the function $B(t,x):= \tilde D(t,x_{12})^{-1}\widebar\Psi(t,x_{11},x_{12})^\top x_2$ with $\tilde D(t,x_{12}) \equiv D(q)$, $x_1:=\col[s,\tq]$, $x_2:=\tilde\theta$ and $\widebar\Psi(t,x_{11},x_{12}) \equiv \widetilde \Psi(q,\dq,\dqr,\ddqr)$. Therefore, since $B_\circ(t,x_2)\equiv D(t,0)^{-1}\widebar\Psi(t,0,0)$ and $D(t,0)$ is positive definite and uniformly bounded for all $t$, the u\ped\ condition on $B_\circ(t,x_2)$ is equivalent to a PE condition on $\widebar\Psi(t,0,0)$ which actually corresponds to $\widetilde \Psi(\qd(t),\dqd(t),\dqd(t),\ddqd(t))$. See \cite{ICRA03} for details.

\section{Conclusions}
\label{sec:conc}
In previous papers we introduced the concept of $\delta$-persistency of excitation as a sufficient condition for uniform convergence of nonlinear systems. In this paper we presented a new mathematical definition which is not formulated with respect to trajectories. Moreover, we have proved for this new definition of  $\delta$-persistency of excitation that it is also necessary for uniform attractivity of certain nonlinear time-varying systems.  Future research includes the use of this property in adaptive control of nonlinearly parameterized systems and persistently excited state observers.

\appendix

\section{Proofs of Lemmas}

\begin{proof}[of Lemma \ref{lem:equiv:defs}]
\nin ({\bf\ref{stmt:B}}) $\Longrightarrow$ u\ped\,: Let $x^*$ be an arbitrary element of ${\cal D}_1$. Define $T^{*}:=T(x^{*})< \infty$ and $\mu^{*} = \mu(x^{*})>0$ so that \rref{eq:B} holds.
Using the uniform continuity of $\phi$, let $\delta>0$
be such that, for all $t \in \mR$,
\begin{equation}
  \norm{z - x^{*}} \leq \delta
   \qquad \Longrightarrow \qquad\norm{\phi(t,z)-\phi(t,x^{*})} \leq  
\frac{\mu^{*}}{2 T^{*}}
    \ .
\end{equation} 
We then have for $z$ such that $\norm{z - x^{*}} \leq \delta$,
\begin{eqnarray}
  \int_{t}^{t+T^{*}} \norm{\phi(\tau,z)} d \tau
    & \geq & \mu^{*} - \int_{t}^{t+T^{*}} \norm{\phi(\tau,z)-\phi(\tau,x^{*})} d \tau 
    \ \geq \ \frac{\mu^{*}}{2} \ .
\end{eqnarray}
\nin u\ped\ $\Longrightarrow$ ({\bf\ref{stmt:B}})\,: Obvious (and uniform continuity is not used).
\end{proof}

\begin{proof}[of Lemma \ref{TA:equiv2:B}]
u\ped\ $\Longrightarrow$ \rref{stmt:A}: Let the u\ped\ property of $\phi$ generate $\delta_x$, $T_x$ and $ \mu_x$. Fix $\Delta\geq\delta>0$ arbitrarily and define $ {\cal D}_{\delta,\Delta}:=\{ x\in\mR^n \, : \, \norm{x_1}\in [\delta,\Delta],\ \norm{x_2}\in [0,\Delta] \}$. Since ${\cal D}_{\delta,\Delta} \subset {\cal D}_1$ we have that \rref{eq:newdpe} holds for any $x\in {\cal D}_{\delta,\Delta}$. Next, we note that the open sets ${\cal B}(x,\delta_x):=\{ z\in\mR^n \, : \, \norm{z-x} < \delta_x\}$ cover ${\cal D}_{\delta,\Delta}$ if we let $x$ range over ${\cal D}_1$ that is, 
\[
{\cal D}_{\delta,\Delta} \subseteq \bigcup_{x\in {\cal D}_1} {\cal B}(x,\delta_{x}) \,.
\]
Since ${\cal D}_{\delta,\Delta}$ is compact from this cover we can extract a finite subcover of ${\cal D}_{\delta,\Delta}$ (see e.g. \cite[ch. 3]{KOLFOM}, \cite[Lemma 5.1, p. 165]{MUN}). In other words, there exists a finite number $n$, and a set $M:=\{x_1, \, x_2, \, \ldots,\, x_n\}$ and corresponding balls ${\cal B}_i(x_i,\delta_{x_i})$ such that the finite union of all these ${\cal B}_i$'s contains ${\cal D}_{\delta,\Delta}$. The statement of the lemma follows defining 
\[
T := \max_{x_i \in M}  T_{x_i} \qquad \mu := \min_{x_i\in M} \mu_{x_i}
\]

\nin  \rref{stmt:A} $\Longrightarrow$ u\ped: Let ({\bf\ref{stmt:A}}) hold with $\bar\Delta$, $\bar\delta>0$, $\bar T$, $\bar \mu$. For any $x\in {\cal D}_{3\bar\delta/2,\Delta}$ we have that \rref{eq:E} holds with $T=\bar T$, $\mu=\bar \mu$ and $\delta=\bar\delta/2$. The result follows by considering $\bar \delta$ arbitrarily small and $\Delta$ arbitrarily large.
\end{proof}

\begin{proof}[of Lemma \ref{lem:equiv:PE}]
By Lemma \ref{TA:equiv2:B} the function $\phi(\cdot,\cdot)$ is u\ped\ if and only if ({\bf\ref{stmt:A}}) holds. We now prove that ({\bf\ref{stmt:A}}) is equivalent to ({\bf\ref{stmt:C}}). 

That ({\bf\ref{stmt:C}}) implies \rref{stmt:A} follows straightforward by defining for each pair $\Delta\geq\delta>0$,\, $T:=\theta_\Delta(\delta)$ and $\mu:=\gamma_\Delta(\delta)$. 

We show now that statement \rref{stmt:A} implies \rref{stmt:C}. Fix $\Delta>0$ otherwise arbitrarily. Let this $\Delta$ and each $\Delta\geq \delta > 0$ generate via Statement \ref{stmt:A}, $\mu_\Delta(\delta)$ and $T_\Delta(\delta)$. Obviously $\mu_\Delta(\cdot)$ and $T_\Delta(\cdot)$ are positive. We define the functions $\bar\mu_\Delta:\mR_{>0}\to\mRp$ and $\widebar T_\Delta:\mR_{>0}\to\mR_{>0}$ as\remfootnote{We define these functions with the lower limit in the inf and sup, being $\min{\delta,\Delta}$ since otherwise $\bar \mu$ and $\bar T$ would have domains only on $(0,\Delta]$. These definitions are the same as saying that for all $\delta>\Delta$ take $\bar{T}(\delta)=\bar{T}(\Delta)$\ldots  }
\begin{equation}\label{infsup}
\bar \mu_\Delta(\delta):= \,\inf_{   \min\{\delta,\Delta\} \leq r \leq \Delta} \mu_\Delta(r) \qquad   \widebar T_\Delta(\delta)  := \sup_{\min\{\delta,\Delta\}  \leq r \leq \Delta} T_\Delta(r)\,.
\end{equation}
 Roughly, $\bar \mu_\Delta(\delta)$ is the smallest over all $\mu_\Delta(\delta)$'s for which \rref{eq:E} holds. Correspondingly, $\widebar T_\Delta$ is the largest of all applicable $T_\Delta$'s. Notice that $\bar \mu_\Delta(\cdot)$ is nondecreasing\remfootnote{since it corresponds to the infimum over shrinking intervals} while $\widebar T_\Delta(\cdot)$ is nonincreasing\remfootnote{since it's a sup over shrinking intervals.} and neither function is necessarily continuous.

Also, from \rref{stmt:A} we have that for each $\delta^*>0$,
\begin{equation}
  \lim\sup_{\shiftleft{6mm}\delta \downarrow \delta^{*}} T_\Delta(\delta) < \infty
\label{eq:14}
\end{equation}
and
\begin{equation}
  \liminf_{\delta \downarrow \delta^{*}} \mu_\Delta(\delta) > 0 \ .
\label{eq:15}
\end{equation}
This follows because if there exists $\delta^*$ such that either \rref{eq:14} or \rref{eq:15} fail to hold then, invoking the definition of u\ped\ with $\delta:= \delta^*/2 $, we obtain that there exists $\beta>0$ and $\theta>0$ such that 
\[
\left\{\, \norm{x_1} \in [\delta^*/2,\,\Delta]\,,\quad \norm{x_2} \in [0,\,\Delta]\,\right\} \ \ \Longrightarrow \ \ \int_t^{t+\theta} \norm{\phi(\tau,x)}d\tau \geq \beta\,.
\]
We conclude that on a neighborhood of $\delta^*$ we can take $\mu_\Delta(\delta):=\min\{\mu_\Delta(\delta^*), \, \beta\} > 0$ and $ T_\Delta(\delta):=\max\{T_\Delta(\delta^*),\,\theta \}$ so that neither \rref{eq:14} nor \rref{eq:15} will fail.

Now we turn the functions $\bar\mu_\Delta(\cdot)$ and $\widebar T_\Delta(\cdot)$ into continuous functions $\bar\gamma_\Delta:\mR_{>0}\to\mRp$ and $\bar\theta_\Delta:\mR_{>0}\to\mR_{>0}$. We define, for each $s > 0$
\begin{equation}
   \bar\gamma_\Delta(s) := \inf_{z > 0}
                        \left\{ \bar\mu_\Delta(z) + |s-z| \right\}
\end{equation}
\begin{equation}
   U_{c}(s) : = \inf_{z > 0}
                          \left\{ \frac{1}{\max\left\{1,\widebar T_\Delta(z)\right\}} + |s-z| \right\} 
\end{equation}
and
\begin{equation}
   \bar\theta_\Delta(s) = \frac{1}{U_{c}(s)} \ .
\end{equation}
If follows from (\ref{eq:14}), (\ref{eq:15}) and the definitions above that $\bar\gamma_\Delta$ and $U_{c}$, and in turn $\bar\theta_\Delta$, take values in $\mR_{>0}$. Invoking standard arguments from nonsmooth analysis we obtain that $\bar\gamma_\Delta$ and $U_{c}$, and in turn $\bar\theta_\Delta$, are continuous (see e.g. \cite[p. 44]{CLARKE}). In addition, from \rref{infsup} we see that $\bar\gamma_\Delta(\cdot)$ is  nondecreasing and  $\bar\theta(\cdot)$ is nonincreasing. 
Moreover,
\begin{equation}\label{nine}
   \bar\gamma_\Delta(s) \leq \bar\mu_\Delta(s)
\end{equation}
and
\begin{equation}
   U_{c}(s) \leq \frac{1}{\max\left\{1,\widebar T_\Delta(s)\right\} }
\end{equation}
so that
\begin{equation}\label{eleven}
  \widebar T_\Delta(s) \leq \max\left\{1,\widebar T_\Delta(s)\right\} \leq  \frac{1}{U_{c}(s)} = \bar\theta_\Delta(s)\,.
\end{equation}
Consequently, in view of \rref{infsup}, \rref{nine} and \rref{eleven}, we have that \rref{eq:E} holds with  $T=\bar\theta_\Delta(\delta)$ and $\mu=\bar\gamma_\Delta(\delta)$. Let  $\gamma_\Delta\in\cK$ and $\theta_\Delta:\mR_{>0}\to\mR_{>0}$ be continuous strictly decreasing such that $\gamma_\Delta(s)\leq \bar\gamma_\Delta(s)$, $\bar\theta_\Delta(s)\leq \theta_\Delta(s)$ for all $s> 0$, $\gamma_\Delta(0)=\bar\gamma_\Delta(0)=0$  and consider any  $x_1\in {\cal B}(\Delta)\backslash \{x_1=0\}$. Then, \rref{eq:E}  holds  with $T=\theta_\Delta(s)$ and $\mu=\gamma_\Delta(s)$, $s=\norm{x_1}$ and $\delta$ arbitrarily small. Statement \rref{stmt:C} follows. 
\end{proof}

\bibliographystyle{plain}
\bibliography{refs}

\newcommand{\SortNoop}[1]{} 
  \def\nesic{Ne\v{s}i\'{c}\,} %
  \def\astrom{{\SortNoop{As}\AA}str{\"{o}}m\,}\let\c=\cedille
\begin{thebibliography}{10}

\bibitem{AEYSEPPEUT}
D.~Aeyels, R.~Spulchre, and J.~Peuteman.
\newblock Asymptotic stability for time-variant systems and observability:
  uniform and nonuniform criteria.
\newblock {\em Math. of Cont. Sign. and Syst.}, 11:1--27, 1998.

\bibitem{CLAN8}
B.~D.~O. Anderson, R.R. Bitmead, C.R. J{ohnson, Jr.}, P.V. Kokotovi{\'c}, R.L.
  Kosut, I.~Mareels, L.~Praly, and B.D. Riedle.
\newblock {\em Stability of adaptive systems}.
\newblock The MIT Press, {Cambridge, MA, USA}, 1986.

\bibitem{AND77}
B.~O. Anderson.
\newblock Exponential stability of linear equations arising in adaptive
  identification.
\newblock {\em IEEE Trans. on Automat. Contr.}, 22(1):83--88, 1977.

\bibitem{ARST78}
Z.~Artstein.
\newblock Uniform asymptotic stability via the limiting equations.
\newblock {\em J. Diff. Eqs.}, 27:172--189, 1978.

\bibitem{ARST82}
Z.~Artstein.
\newblock Stability, observability and invariance.
\newblock {\em J. Diff. Eqs.}, 44:224--248, 1982.

\bibitem{ASTBOH}
K.~J. \astrom and Bohn.
\newblock Numerical identification of linear dynamic systems from normal
  operating records.
\newblock In P.~H. Hammond, editor, {\em {{\em Proc. of the 2nd IFAC Symp. on
  Theory of Self-adaptive Control Systems}}}, pages 96--111, {Nat. Phys. Lab.,
  Teddington, England}, 1965.

\bibitem{BRO}
R.~Brockett.
\newblock Asymptotic stability and feedback stabilization.
\newblock In R.~S.~Millman R.~W.~Brocket and H.~J. Sussmann, editors, {\em
  Differential geometric control theory}, pages 181--191. {Birkh\"auser}, 1983.

\bibitem{CLARKE}
F.~H. Clarke, Yu.~S. Ledyaev, R.~J. Stern, and P.~R. Wolenski.
\newblock {\em Nonsmooth analysis and control theory}.
\newblock {G}raduate {T}exts in {M}athematics. Springer-Verlag, 1998.

\bibitem{THORSTUFFIJRNLC}
T~I. Fossen, A.~\loria\, and A.~Teel.
\newblock A theorem for {UGAS} and {ULES} of (passive) nonautonomous systems:
  Robust control of mechanical systems and ships.
\newblock {\em Int. J. Rob. Nonl. Contr.}, 11:95--108, 2001.

\bibitem{THOROBS}
T.~I. Fossen and J.~P. Strand.
\newblock Passive nonlinear observer design for ships using {Lyapunov} methods:
  Full-scale experiments with a supply vessel.
\newblock {\em Automatica}, 35(1):3--16, 1999.

\bibitem{GAUPOI}
M.~Gauthier and P.~Poignet.
\newblock Identification nonlin\'eaire continue en boucle ferm\'ee des
  param\`etres physiques de syst\`emes m\'ecatroniques par mod\`ele inverse et
  moindres carr\'es d'erreur d'entr'e.
\newblock In {\em Journ\'ees d'identification et mod\'elisation
  exp\'erimentale}, {Vandoeuvre-l\`es-Nancy, France}, Mars 2001.

\bibitem{HALE}
J.K. Hale.
\newblock {\em Ordinary Differential equations}.
\newblock Interscience. {John Wiley}, {New York}, 1969.

\bibitem{JANIJACSP}
M.~Jankovi\'c.
\newblock Adaptive output feedback control of nonlinear feedback linearizable
  systems.
\newblock {\em Int. J. Adapt. Contr. Sign. Process.}, 10(1):1--18, 1996.

\bibitem{JIANONHOL96}
Z.~P. Jiang.
\newblock Iterative design of time-varying stabilizers for multi-input systems
  in chained form.
\newblock {\em Syst. \& Contr. Letters}, 28:255--262, 1996.

\bibitem{KALMAN}
R.~E. Kalman.
\newblock Contributions to the theory of optimal control.
\newblock {\em Bol. Soc. Mat. Mexicana}, 5:102--119, 1960.

\bibitem{KHALIL}
H.~Khalil.
\newblock {\em Nonlinear systems}.
\newblock {Macmillan Publishing Co., 2nd ed.}, New York, 1996.

\bibitem{KOLFOM}
A.~N. Kolmogorov and S.~V. Fomin.
\newblock {\em Introductory real analysis}.
\newblock Dover, {Mineola, N.Y.}, 1970.
\newblock {ISBN:} 0-486-61226-0.

\bibitem{TCLEETAC01}
T.~C. Lee and B.~S. Chen.
\newblock A general stability criterion for time-varying systems using a
  modified detectability condition.
\newblock {\em IEEE Trans. on Automat. Contr.}, 47(5):797--802, 2002.

\bibitem{ANA99}
A.~P. Loh, A.~M. Annaswamy, and F.~P. Skantze.
\newblock Adaptation in the prescence of a general nonlinear parameterization:
  an error model approach.
\newblock {\em IEEE Trans. on Automat. Contr.}, 44(9):1634--1653, 1999.

\bibitem{ICRA03}
A.~Loria, R.~Kelly, and A.~Teel.
\newblock Uniform parametric convergence in the adaptive control of
  manipulators: a case restudied.
\newblock In {\em Proc. IEEE Conf. Robotics Automat.}, {Taipei, Taiwan}, 2003.
\newblock Submitted.

\bibitem{NONHOLEJC}
A.~\loria, E.~Panteley, and K.~Melhem.
\newblock {UGAS} of ``skew-symmetric'' time-varying systems: application to
  stabilization of chained form systems.
\newblock {\em EJC}, 8(1):33--43, 2002.

\bibitem{MATCDC02}
A.~\loria\, E.~Panteley, D.~Popovi\'c, and A.~Teel.
\newblock An extension of matrosov's theorem with application to nonholonomic
  control systems.
\newblock In {\em Proc. 40th. IEEE Conf. Decision Contr.}, {Las Vegas, CA,
  USA}, 2002.
\newblock {Paper no. REG0625}.

\bibitem{uped}
A.~\loria\, E.~Panteley, and A.~Teel.
\newblock A new persistency-of-excitation condition for {UGAS} of {NLTV}
  systems: {Application} to stabilization of nonholonomic systems.
\newblock In {\em Proc. 5th. European Contr. Conf.}, 1999.
\newblock Paper no. 500.

\bibitem{MARTOM}
R.~Marino and P.~Tomei.
\newblock Global adaptive output feedback control of nonlinear systems. {Part
  I~}: {L}inear parameterization.
\newblock {\em IEEE Trans. on Automat. Contr.}, 38:17--32, 1993.

\bibitem{MAZPRA}
F.~Mazenc and L.~Praly.
\newblock Adding integrators, saturated controls and global asymptotic
  stabilization of feedforward systems.
\newblock {\em IEEE Trans. on Automat. Contr.}, 41(11):1559--1579, 1996.

\bibitem{MAZPRA99}
F.~Mazenc and L.~Praly.
\newblock Asymptotic tracking of a reference state for systems with a
  feedforward structure.
\newblock {\em Automatica}, 36:179--187, 1999.

\bibitem{MORNAR}
A.~P. Morgan and K.~S. Narendra.
\newblock {On the stability of nonautonomous differential equations $\dot x =
  [A+B(t)]x$ with skew-symmetric matrix $B(t)$}.
\newblock {\em SIAM J. on Contr. and Opt.}, 15(1):163--176, 1977.

\bibitem{MORNAR2}
A.~P. Morgan and K.~S. Narendra.
\newblock On the uniform asymptotic stability of certain linear nonautonomous
  differential equations.
\newblock {\em SIAM J. on Contr. and Opt.}, 15(1):5--24, 1977.

\bibitem{MORSAMEJC}
P.~Morin and C.~Samson.
\newblock Application of backstepping techniques to the time-varying
  exponential stabilization of chained-form systems.
\newblock {\em European J. of Contr.}, 3(1):15--37, 1997.

\bibitem{SAMSIAM02}
P.~Morin and C.~Samson.
\newblock A characterization of the {Lie} algebra rank condition by transverse
  periodic functions.
\newblock {\em SIAM J. on Contr. and Opt.}, 40(4):1227--1249, 2002.

\bibitem{MUN}
J.~R. Munkres.
\newblock {\em Topology: a first course}.
\newblock {Prentice-Hall, Inc.,}, {Englewood Cliffs, NJ}, 1975.

\bibitem{MURSAS}
R.~M. Murray and S.~S. Sastry.
\newblock {Nonholonomic motion planning: Steering with sinusoids}.
\newblock {\em IEEE Trans. on Automat. Contr.}, 38(5):700--716, 1993.

\bibitem{NARANA}
K.~S. Narendra and A.~M. Annaswamy.
\newblock {\em Stable adaptive systems}.
\newblock {Prentice-Hall, Inc.}, New Jersey, 1989.

\bibitem{OrtegaFradkov}
R.~Ortega and A.~L. Fradkov.
\newblock Asymptotic stability of a class of adaptive systems.
\newblock {\em Int. J. Adapt. Contr. Sign. Process.}, 7:255--260, 1993.

\bibitem{PBCELS}
R.~Ortega, A.~\loria\, P.~J. Nicklasson, and H.~Sira-Ram\'{\i}rez.
\newblock {\em {Passivity-based Control of Euler-Lagrange Systems: Mechanical,
  Electrical and Electromechanical Applications}}.
\newblock Comunications and Control Engineering. Springer Verlag, London, 1998.
\newblock ISBN 1-85233-016-3.

\bibitem{ORTSPO}
R.~Ortega and M.~Spong.
\newblock Adaptive motion control of rigid robots: A tutorial.
\newblock {\em Automatica}, 25-6:877--888, 1989.

\bibitem{TACDELTAPE}
E.~Panteley, A.~\loria, and A.~Teel.
\newblock Relaxed persistency of excitation for uniform asymptotic stability.
\newblock {\em IEEE Trans. on Automat. Contr.}, 46(12):1874--1886, 2001.

\bibitem{POM92}
J.~B. Pomet.
\newblock Explicit design of time-varying stabilizing control laws for a class
  of controllable systems without drift.
\newblock {\em Syst. \& Contr. Letters}, 18:147--158, 1992.

\bibitem{ROUMAW}
N.~Rouche and J.~Mawhin.
\newblock {\em {Ordinary differential equations II: Stability and periodical
  solutions}}.
\newblock Pitman publishing Ltd., London, 1980.

\bibitem{SAM1st}
C.~Samson.
\newblock Time-varying stabilization of a car-like mobile robot.
\newblock Technical report, {INRIA Sophia-Antipolis}, 1990.

\bibitem{SAM3}
C.~Samson.
\newblock {Control of chained system: Application to path following and
  time-varying point stabilization of mobile robots}.
\newblock {\em IEEE Trans. on Automat. Contr.}, 40(1):64--77, 1995.

\bibitem{SASBOD}
S.~Sastry and M.~Bodson.
\newblock {\em Adaptive control: Stability, convergence and robustness}.
\newblock {Prentice Hall Intl.}, 1989.

\bibitem{SCISIC}
L.~Sciavicco and B.~Siciliano.
\newblock {\em Modeling and control of robot manipulators}.
\newblock {McGraw Hill}, {New York}, 1996.

\bibitem{SEPBOOK}
R.~Sepulchre, M.~Jankovi\'c, and P.~Kokotovi\'c.
\newblock {\em Constructive nonlinear control}.
\newblock {Springer Verlag}, 1997.

\bibitem{SLOTLIA}
J.~J. Slotine and W.~Li.
\newblock Adaptive manipulator control: a case study.
\newblock {\em IEEE Trans. on Automat. Contr.}, AC-33:995--1003, 1988.

\bibitem{SLOTLIPE}
J.J. Slotine and W.~Li.
\newblock Theoretical issues in adaptive manipulator control.
\newblock In {\em {5th Yale Workshop on Apl. Adaptive Systems Theory}}, pages
  252--258, 1987.

\bibitem{SPOVID}
M.~Spong and M.~Vidyasagar.
\newblock {\em Robot Dynamics and Control}.
\newblock John Wiley \& Sons, New York, 1989.

\bibitem{SPOORTKEL}
M.~W. Spong, R.~Ortega, and R.~Kelly.
\newblock {Comments on "Adaptive Manipulator Control: A Case Study".}
\newblock {\em IEEE Trans. on Automat. Contr.}, 35(6):761, 1990.

\bibitem{INTLEMMCSS}
A.~Teel, E.~Panteley, and A.~\loria.
\newblock Integral characterizations of uniform asymptotic and exponential
  stability with applications.
\newblock {\em Math. of Cont. Sign. and Syst.}, 15:177--201, 2002.

\bibitem{VORIFAC02}
V.~I. Voritkonov.
\newblock Partial stability, stabilization and control: a some recent results.
\newblock In {\em Proc. 15th. IFAC World Congress}, {Barcelona, Spain}, 2002.
\newblock CDROM file: {\tt 02442.pdf}.

\end{thebibliography}


\end{document}